\documentclass[12pt,a4paper,fleqn]{article}
\usepackage{a4wide,amsfonts,amsmath,latexsym,amssymb,euscript,graphicx,units,mathrsfs}

\usepackage{graphicx}
\usepackage{color}
\usepackage{amssymb}
\usepackage{amssymb}
\usepackage[T1]{fontenc}
\usepackage{latexsym}
\usepackage{xypic}
\usepackage{eufrak}
\usepackage{euscript}
\usepackage{amsfonts,amsmath}
\usepackage{verbatim}
\usepackage{fancyhdr}
\usepackage[english]{babel}
\usepackage{mathrsfs}
\usepackage{units}

\newtheorem{prop}{Proposition}[section]
\newtheorem{cor}[prop]{Corollary}
\newtheorem{lemme}[prop]{Lemma}
\newtheorem{rem}[prop]{Remark}
\newtheorem{rems}[prop]{Remarks}
\newtheorem{thm}[prop]{Theorem}
\newtheorem{defi}[prop]{Definition}

\newtheorem{probl}[prop]{Problem}

\renewcommand{\geq}{\geqslant}
\def\leq{\leqslant}

\newcommand{\Z}{\mathbb{Z}}
\newcommand{\R}{\mathbb{R}}

\newcommand{\T}{\mathbb{T}}

\def\HH{\EuFrak H}

\font\tenbb=msbm10

\def\1{{\mathbf{1}}}

\font\tenbb=msbm10

\def\zZ{\hbox{\tenbb Z}}

\def\1{{\mathbf{1}}}
\def\0.5{{\frac{1}{2}}}

\newcommand{\fin}
{ \vspace{-0.6cm}
\begin{flushright}
\mbox{$\Box$}
\end{flushright}
\noindent }

\newcommand{\qed}{\nopagebreak\hspace*{\fill}
{\vrule width6pt height6ptdepth0pt}\par}

\newcounter{rea}
\begin{document}

\begin{center}
{\Large{\bf Optimal Berry-Esseen rates on the Wiener space: the barrier of third and fourth cumulants}}
\normalsize
\\~\\ by Hermine Bierm\'e \footnote{Email: {\tt hermine.bierme@mi.parisdescartes.fr}; HB was partially supported by the ANR grant
`MATAIM' ANR-09-BLAN-0029-01}, 
Aline Bonami \footnote{Email: {\tt Aline.Bonami@univ-orleans.fr}; AB was partially supported by the ANR grant
`MATAIM' ANR-09-BLAN-0029-01}, 
Ivan Nourdin\footnote{Email: {\tt inourdin@gmail.com}; IN was partially supported by the 
ANR grants ANR-09-BLAN-0114 and ANR-10-BLAN-0121.},
Giovanni Peccati\footnote{Email: {\tt giovanni.peccati@gmail.com}}\\ 
{\it Universit\'e Paris Descartes and Universit\'e de Tours, Universit\'e d'Orl\'eans, Universit\'e de Nancy and Fondation Sciences Math\'ematiques de Paris, Universit\'e du Luxembourg}\\~\\
\end{center}

{\small \noindent {{\bf Abstract}: Let $\{F_n : n\geq 1 \}$ be a normalized sequence of random variables in some fixed Wiener chaos associated with
a general Gaussian field, and assume that $E[F_n^4] \rightarrow E[N^4]=3$, where $N$ is a standard Gaussian random variable. Our main result is the following general bound: there exist two finite constants $c,C>0$ such that,
for $n$ sufficiently large, $c\times \max(|E[F_n^3]|, E[F_n^4]-3)  \leq d(F_n,N) \leq C\times  \max(|E[F_n^3]|, E[F_n^4]-3) $, where
$d(F_n,N) = \sup|E[h(F_n)] - E[h(N)]|$, and $h$ runs over the class of all real functions with a second derivative bounded by 1. This
shows that the deterministic sequence $\max(|E[F_n^3]|, E[F_n^4]-3)$, $n\geq 1$, completely characterizes the rate of convergence (with
respect to smooth distances) in CLTs involving chaotic random variables. These results are used to determine optimal rates of convergence in the 
Breuer-Major central limit theorem, with specific emphasis on fractional Gaussian noise.
\\



\section{Introduction and main results}

Let $X = \{X(h) : h\in \HH\}$ be an isonormal Gaussian process (defined on an adequate space $(\Omega, \mathcal{F}, P)$) over some real separable Hilbert space $\HH$, fix an integer $q\geq 2$, and let $\{F_n : n\geq 1\}$ be a sequence of random variables belonging to the $q$th Wiener chaos of $X$ (see Section \ref{ss:isonormal} for precise definitions). Assume that $E[F_n^2] = 1$ for every $n$. In recent years, many efforts have been devoted to the characterization of those chaotic sequences $\{F_n\}$ verifying a Central Limit Theorem (CLT), that is, such that $F_n$ converges in distribution to $N\sim \mathscr{N}(0,1)$ (as $n\to \infty$), where $\mathscr{N}(0,1)$ denotes a centered Gaussian law with unit variance. An exhaustive solution to this problem was first given by Nualart and Peccati in \cite{nunugio}, in the form of the following ``fourth moment theorem''.

\begin{thm}[Fourth Moment Theorem -- see \cite{nunugio}]\label{t:nunugio} Fix an integer $q\geq 2$, and consider a sequence of random variables
$\{F_n : n\geq 1\}$ belonging to the $q$th Wiener chaos of $X$ and such that $E[F_n^2]=1$
for all $n\geq 1$. Then, as $n\rightarrow \infty$, $F_n$ converges in distribution to $N\sim\mathscr{N}(0,1)$ if and only if $E[F_n^4] \rightarrow E[N^4] = 3$ .
\end{thm}

Note that Theorem \ref{t:nunugio} represents a drastic simplification of the usual {\it method of moments and cumulants}, as described e.g. in \cite{PeTa}. Combining the so-called {\it Stein's method} for normal approximations (see \cite{ChenGoldShao, np-book}, as well as Section \ref{ss:stein} below) with Malliavin calculus (see \cite{MallBook, nualartbook}, as well as Section \ref{ss:mall}), one can also prove the forthcoming Theorem \ref{t:npbounds}, providing explicit upper bounds in the {\it total variation distance}. We recall that the total variation distance $d_{TV}(F,G)$ between the laws of two real-valued random variables $F,G$ is defined as
\[
d_{TV}(F,G) = \sup_{A\in \mathscr{B}(\R)} \left| P[F\in A] -P[G\in A]\right|,
\]
where the supremum runs over the class of all Borel sets $A\subset \R$. Given a smooth functional of the isonormal process $X$, we shall also write $DF$ to indicate the Malliavin derivative of $F$ (thus $DF$ is a $\HH$-valued random element -- see again Section \ref{ss:mall} for details).
\begin{thm}[Fourth Moment Bounds -- see \cite{np-ptrf, npr-aop}]\label{t:npbounds} Fix $q\geq 2$, let $F$ be an element of the $q$th Wiener chaos of $X$ with unit variance, and let $N\sim \mathscr{N}(0,1)$. The following bounds are in order:
\begin{equation}\label{e:npupper}
d_{TV}(F,N) \leq 2 \sqrt{E\left[ \left(1 - \frac1q\|DF\|^2_\HH \right)^2\right] } \leq 2\sqrt{\frac{q-1}{3q} \left(E[F^4]-3\right) } .
\end{equation}

\end{thm}

\begin{rem}{\rm
\begin{enumerate}

\item The two inequalities in (\ref{e:npupper}) were discovered, respectively, in \cite{np-ptrf} and \cite{npr-aop}. Using the properties of the Malliavin derivative $DF$ (see Section \ref{ss:mall} below), one sees immediately that
\[
E\left[ \left(1 - \frac1q\|DF\|^2_\HH \right)^2\right] = {\rm Var}\left(\frac1q\|DF\|^2_\HH\right).
\]
\item One can prove the following refinement of the second inequality in (\ref{e:npupper}) (see \cite[Lemma 3.5]{survey}): for every random variable $F$ belonging to the $q$th Wiener chaos of $X$ and with unit variance
\begin{equation}\label{e:vev}
{\rm Var}\left(\frac1q\|DF\|^2_\HH\right)
\leq \frac{q-1}{3q}\big(E[F^4]-3\big)
\leq (q-1){\rm Var}\left(\frac1q\|DF\|^2_\HH\right).
\end{equation}

\item Theorem \ref{t:npbounds} implies that, not only the condition $E[F_n^4]\to 3$ is necessary and sufficient for convergence to Gaussian, as stated in Theorem \ref{t:nunugio}, but also that the sequences
\begin{equation}\label{e:betagamma}
\beta(n) := \sqrt{E[F_n^4] - 3}\,\,\,\, \mbox{\ and \ } \gamma(n) := \sqrt{{\rm Var}\left(\frac1q\|DF_n\|^2_\HH\right)}, \quad n\geq 1,
\end{equation}
bound from above (up to a constant) the speed of convergence of the law of $F_n$ to that of $N$ in the topology induced by $d_{TV}$.
\item If one replaces the total variation distance with the Kolmogorov distance or with the Wasserstein distance (see e.g. \cite{np-ptrf, np-book} for definitions), then the bounds (\ref{e:npupper}) hold without the multiplicative factor $2$ before the square roots.
\item When $E[F] = 0$ and $E[F^2]=1$, the quantity $E[F^4] - 3$ coincides with the {\it fourth cumulant} of $F$, see Definition \ref{D : cum}. One can also prove that, if $F$ belongs to a fixed Wiener chaos and has unit variance, then $E[F^4] > 3$ (see \cite{nunugio}). 
\item Throughout the paper, in order to simplify the notation, we only consider sequences of random variables having unit variance. The extension to arbitrary sequences whose variances converge to a constant can be deduced by a straightforward adaptation of our arguments.
\end{enumerate}
}
\end{rem}

A natural problem is now the following. 

\begin{probl}\label{p:uno}
{Assume that $\{F_n\}$ is a unit variance sequence belonging to the $q$th Wiener chaos of the isonormal Gaussian process $X$. Suppose that $F_n$ converges in distribution to $N\sim \mathscr{N}(0,1)$ and fix a distance $d_0(\cdot,\cdot)$ between the laws of real-valued random variables. Can one find an explicit optimal rate of convergence associated with the distance $d_0$?} 
\end{probl}

The notion of optimality adopted throughout the paper is contained in the next definition.

\begin{defi}\label{d:optimal}{\rm
Assume that, as $n\to\infty$, $F_n$ converges in distribution to $N$, and fix a generic distance $d_0(\cdot,\cdot)$ between the laws of 
real-valued random variables. A deterministic sequence $\{\varphi(n) : n\geq 1\}$ such that $\varphi(n)\downarrow 0$ is said to provide an 
{\it optimal rate of convergence} with respect to $d_0$ if there exist constants $0<c<C<\infty$ (not depending on $n$) 
such that, for $n$ large enough,
\begin{equation}\label{e:optr}
c \leq \frac{d_{0}(F_n,N)}{\varphi(n)} \leq C. 
\end{equation}
}
\end{defi}

The problem of finding optimal rates is indeed quite subtle. A partial solution to Problem \ref{p:uno} is given by Nourdin and Peccati in \cite{np-aop}. In this reference, a set of sufficient conditions are derived, ensuring that the sequences $\beta(n),\gamma(n)$ in (\ref{e:betagamma}) yield optimal rates for the distance $d_0 = d_{TV}$.  In particular, these conditions involve the joint convergence of the two-dimensional vectors 
\begin{equation}\label{e:2dimv}
\left(F_n, \frac{1-q^{-1}\|DF_n\|^2_\HH}{\gamma(n)}\right),\quad n\geq 1. 
\end{equation}
The following statement constitutes one of the main finding of \cite{np-aop} (note that the reference \cite{np-aop} only deals with the Kolmogorov distance but, as far as lower bounds are concerned, it is no more difficult to work directly with $d_{TV}$).

\begin{thm}[See \cite{np-aop}]\label{t:npoptimal}
Let $\{F_n\}$ be a unit variance sequence belonging to the $q$th Wiener chaos of $X$, and suppose that, as $n\to\infty$, $F_n$ converges in distribution to $N\sim \mathscr{N}(0,1)$. Assume moreover that the sequence of two-dimensional random vectors in (\ref{e:2dimv}) converges in distribution to a Gaussian vector $(N_1,N_2)$ such that $E[N_1^2] =E[N_2^2]=1$, and $E[N_1N_2] =:\rho \in (-1,1)$. 
Then, for every $z\in\R$:
\begin{equation}\label{Uber}
\gamma(n)^{-1}[P(F_n\leq z) - P(N\leq z)]
\to
\frac{\rho}{3}(z^2-1)\frac{e^{-z^2/2}}{\sqrt{2\pi}}\quad\mbox{as $n\to\infty$}.
\end{equation}
In particular, if $\rho\neq 0$ the sequences $\beta(n),\gamma(n)$ defined in (\ref{e:betagamma}) provide optimal rates of convergence with respect to the total variation distance $d_{TV}$, in the sense of Definition \ref{d:optimal}.
\end{thm}

As shown in \cite[Theorem 3.1]{np-aop} the conditions stated in Theorem \ref{t:npoptimal} can be generalized to arbitrary sequences of smooth random variables (not necessarily belonging to a fixed Wiener chaos). Moreover, the content of Theorem \ref{t:npoptimal} can be restated in terms of {contractions} (see \cite[Theorem 3.6]{np-aop}) or, for elements of the second Wiener chaos of $X$, in terms of cumulants (see \cite[Proposition 3.8]{np-aop}). 

\smallskip

One should note that the techniques developed in \cite{np-aop} also imply analogous results for the Kolmogorov and the Wasserstein distances, that we shall denote respectively by $d_{Kol}$ and $d_W$. However, although quite flexible and far-reaching, the findings of \cite{np-aop} {\sl do not allow} to deduce a complete solution (that is, a solution valid for arbitrary sequences $\{F_n\}$ in a fixed Wiener chaos) of Problem \ref{p:uno} for either one of the distances $d_{TV}$, $d_{Kol}$ and $d_W$. 
For instance, the results of \cite{np-aop} provide optimal rates in the Breuer-Major CLT only when the involved subordinating 
function has an {\it even} Hermite rank, whereas the general case remained an open problem till now  -- see \cite[Section 6]{np-aop}.

\smallskip

The aim of this paper is to provide an exhaustive solution to Problem \ref{p:uno} in the case of a suitable smooth distance between laws of real-valued 
random variables. The distance we are interested in is denoted by $d(\cdot,\cdot)$, and involves test functions that are twice differentiable. The formal definition of $d$ is given below.

\begin{defi}\label{d:distance}
{\rm Given two random variables $F,G$ with finite second moments we write $d(F,G)$ in order to indicate the quantity
\[
d(F,G)=\sup_{h\in\mathcal{U}}\big|E[h(F)]-E[h(G)]\big|,
\]
where $\mathcal{U}$ stands for the set of functions $h:\R\to \R$ which are
$\mathcal{C}^2$ (that is, twice differentiable and with continuous derivatives) and such that $\|h''\|_\infty\leq 1$. 
}
\end{defi}

Observe that
$d(\cdot,\cdot)$ defines an actual distance on the class of the distributions of random variables having a finite second moment. Also,
the topology induced by $d$ on this class is stronger than the topology of the convergence in distribution, that is: if $d(F_n,G)\to 0$,
then $F_n$ converges in distribution to $G$.

The forthcoming Theorem \ref{t:mainupper} and Theorem \ref{t:mainlower} contain the principal upper and lower bounds proved in this work: 
once merged, they show that the sequence
\begin{equation}\label{maxbig}
\max\big\{|E[F_n^3]|,E[F_n^4]-3\big\}, \quad n\geq 1,
\end{equation}
always provides optimal rates for the distance $d$, whenever $\{F_n\}$ lives inside a fixed Wiener chaos. 
As anticipated, this yields an exhaustive solution to Problem \ref{p:uno} in the case $d_0=d$. 
One should also note that the speed of convergence to zero of the quantity (\ref{maxbig})
can be given by either one of the two sequences $\{| E[F_n^3]|\}$ and $\{E[F_n^4]-3\}$;
see indeed Corollary \ref{cor-surprising} for explicit examples of both situations.

\begin{rem}\label{r:convm}{\rm Let $\{F_n : n\geq 1\}$ be a sequence of random variables living inside a finite sums of Wiener chaoses. Assume that $E[F_n^2]\to 1$ and that $F_n$ converges in distribution to $N\sim\mathscr{N}(0,1)$. Then, the hypercontractivity property of the Wiener chaos (see e.g. \cite[Chapter V]{Janson}) imply that $E[F_n^k] \to E[N^k]$ for every integer $k\geq 3$. In particular, one has necessarily that $E[F_n^3]\to 0$.
}
\end{rem}

\begin{thm}[Upper bounds]\label{t:mainupper}
Let $N\sim\mathscr{N}(0,1)$ be a standard Gaussian random variable.
Then, there exists $C>0$ such that, for all integer $q\geq 2$ and
all element $F$ of the $q$th Wiener chaos with unit variance,
\begin{equation}\label{cl-thm1}
d(F,N)\leq C\max\big\{|E[F^3]|,E[F^4]-3\big\}.
\end{equation}
\end{thm}

\begin{rem}{\rm
\begin{enumerate}
\item In the statement of Theorem \ref{t:mainupper}, we may assume without loss of generality that $N$ is stochastically independent of $F$. 
Then, by suitably using integration by parts and then Cauchy-Schwarz (see e.g. \cite[Theorem 3.2]{NPP}), one can show that, 
for every $h\in\mathcal{U}$ (see Definition \ref{d:distance}),
\begin{eqnarray}
\big|E[h(F)]-E[h(N)]\big|&=&\frac12\left|\int_0^1 E\left[h''(\sqrt{1-t}F+\sqrt{t}N)\left(1-\frac1q\|DF\|^2_\HH\right)\right]dt\right|\notag\\
&\leq&\frac12 \sqrt{{\rm Var}\left(\frac1q\|DF\|^2_\HH\right)}\label{naive}.
\end{eqnarray}
Since (\ref{e:vev}) is in order, one sees that the inequality (\ref{naive}) does not allow to obtain a better bound than
\[
d(F,N)\leq C\sqrt{E[F^4]-3},
\]
which is not sharp in general, compare indeed with (\ref{cl-thm1}). One should observe that the rate $\sqrt{E[F_n^4]-3}$ may happen to be optimal in some instances, precisely when
$E[F_n^3]$ and $\sqrt{E[F_n^4]-3}$ have the same order. In the already quoted paper \cite{np-aop} one can find several explicit examples where this phenomenon takes place.
\item Let $F$ be an element of the $q$th Wiener chaos of some isonormal Gaussian process, and assume that $F$ has variance 1.
It is shown in \cite[Proposition 3.14]{survey} that there exists a
constant $C$, depending only on $q$, such that $|E[F^3]| \leq C \sqrt{E[F^4] - 3}$.
Using this fact, one can therefore obtain yet another proof of Theorem \ref{t:nunugio} based on the upper bound (\ref{cl-thm1}).
\end{enumerate}
}
\end{rem}

\begin{thm}[Lower bounds]\label{t:mainlower}
Fix an integer $q\geq 2$ and consider a sequence of random variables $\{F_n : n\geq 1\}$ belonging to the $q$th Wiener chaos of some
isonormal Gaussian process and such that $E[F_n^2]=1$. Assume that, as $n\to \infty$, $F_n$ converges in distribution to
$N\sim\mathscr{N}(0,1)$. Then there exists 
 $c>0$ (depending on the sequence $\{F_n\}$, but not on $n$) such that
\begin{equation}\label{cl-thm2}
d(F_n,N)\geq c\times\max\big\{|E[F_n^3]|,E[F_n^4]-3\big\}, \quad n\geq 1.
\end{equation}
\end{thm}

Our proofs revolve around several new estimates (detailed in Section \ref{s:estimates}), that are in turn based on the analytic characterization of cumulants given in \cite{np-jfa}. Also, a fundamental role is played by the Edgeworth-type expansions introduced by Barbour in \cite{BarbourPtrf1986}.

\subsection{Plan}

The paper is organized as follows. Section \ref{s:malliavin} is devoted to some preliminary results of Gaussian analysis and Malliavin calculus. Section \ref{s:steincum} deals with Stein's method, cumulants and Edgeworth-type expansions. Section \ref{s:estimates} contains the main technical estimates of the paper. Section \ref{s:proofs} focuses on the proofs of our main findings, whereas in Section \ref{s:BM} one can find several applications to the computation of optimal rates in the Breuer-Major CLT.

\section{Elements of Gaussian analysis and Malliavin calculus}\label{s:malliavin}

This section contains the essential elements of Gaussian analysis and Malliavin calculus that are used in this paper. See the classical
references \cite{MallBook, nualartbook} for further details.

\subsection{Isonormal processes and multiple integrals}\label{ss:isonormal}

Let $\EuFrak H$ be a real separable Hilbert space. For any $q\geq 1$, we write $\EuFrak H^{\otimes q}$ and $\EuFrak H^{\odot q}$ to
indicate, respectively, the $q$th tensor power and the $q$th symmetric tensor power of $\EuFrak H$; we also set by convention
$\EuFrak H^{\otimes 0} = \EuFrak H^{\odot 0} =\R$. When $\HH = L^2(A,\mathcal{A}, \mu) =:L^2(\mu)$, where $\mu$ is a $\sigma$-finite
and non-atomic measure on the measurable space $(A,\mathcal{A})$, then $\EuFrak H^{\otimes q} = L^2(A^q,\mathcal{A}^q,\mu^q)=:L^2(\mu^q)$, and $\EuFrak H^{\odot q} = L_s^2(A^q,\mathcal{A}^q,\mu^q) := L_s^2(\mu^q)$, where $L_s^2(\mu^q)$ stands for the subspace of $L^2(\mu^q)$ composed of those functions that are $\mu^q$-almost everywhere symmetric. We denote by $X=\{X(h) : h\in \EuFrak H\}$
an {\it isonormal Gaussian process} over
$\EuFrak H$.
This means that $X$ is a centered Gaussian family, defined on some probability space $(\Omega ,\mathcal{F},P)$, with a covariance structure given by the relation
$E\left[ X(h)X(g)\right] =\langle h,g\rangle _{\EuFrak H}$. We also assume that $\mathcal{F}=\sigma(X)$, that is,
$\mathcal{F}$ is generated by $X$.

For every $q\geq 1$, the symbol $\mathcal{H}_{q}$ stands for the $q$th {\it Wiener chaos} of $X$,
defined as the closed linear subspace of $L^2(\Omega ,\mathcal{F},P) =:L^2(\Omega) $
generated by the family $\{H_{q}(X(h)) : h\in \EuFrak H,\left\|
h\right\| _{\EuFrak H}=1\}$, where $H_{q}$ is the $q$th Hermite polynomial
given by 
\begin{equation}\label{hq}
H_q(x) = (-1)^q e^{\frac{x^2}{2}}
 \frac{d^q}{dx^q} \big( e^{-\frac{x^2}{2}} \big).
\end{equation}
We write by convention $\mathcal{H}_{0} = \mathbb{R}$. For
any $q\geq 1$, the mapping $I_{q}(h^{\otimes q})=H_{q}(X(h))$ can be extended to a
linear isometry between the symmetric tensor product $\EuFrak H^{\odot q}$
(equipped with the modified norm $\sqrt{q!}\left\| \cdot \right\| _{\EuFrak H^{\otimes q}}$)
and the $q$th Wiener chaos $\mathcal{H}_{q}$. For $q=0$, we write $I_{0}(c)=c$, $c\in\mathbb{R}$. A crucial fact is that, when $\HH = L^2(\mu)$, for every $f\in\EuFrak H^{\odot q} = L_s^2(\mu^q)$ the random variable $I_q(f)$ coincides with the $q$-fold multiple Wiener-It\^o stochastic integral of $f$ with respect to the centered Gaussian measure (with control $\mu$) canonically generated by $X$ (see \cite[Section 1.1.2]{nualartbook}).

It is
well-known that $L^2(\Omega)$
can be decomposed into the infinite orthogonal sum of the spaces $\mathcal{H}_{q}$. It follows that any square-integrable random variable
$F\in L^2(\Omega)$ admits the following {\it Wiener-It\^{o} chaotic expansion}
\begin{equation}
F=\sum_{q=0}^{\infty }I_{q}(f_{q}),  \label{E}
\end{equation}
where $f_{0}=E[F]$, and the $f_{q}\in \EuFrak H^{\odot q}$, $q\geq 1$, are
uniquely determined by $F$. For every $q\geq 0$, we denote by $J_{q}$ the
orthogonal projection operator on the $q$th Wiener chaos. In particular, if
$F\in L^2(\Omega)$ is as in (\ref{E}), then
$J_{q}F=I_{q}(f_{q})$ for every $q\geq 0$.

Let $\{e_{k},\,k\geq 1\}$ be a complete orthonormal system in $\EuFrak H$.
Given $f\in \EuFrak H^{\odot p}$ and $g\in \EuFrak H^{\odot q}$, for every
$r=0,\ldots ,p\wedge q$, the \textit{contraction} of $f$ and $g$ of order $r$
is the element of $\EuFrak H^{\otimes (p+q-2r)}$ defined by
\begin{equation}
f\otimes _{r}g=\sum_{i_{1},\ldots ,i_{r}=1}^{\infty }\langle
f,e_{i_{1}}\otimes \ldots \otimes e_{i_{r}}\rangle _{\EuFrak H^{\otimes
r}}\otimes \langle g,e_{i_{1}}\otimes \ldots \otimes e_{i_{r}}
\rangle_{\EuFrak H^{\otimes r}}.  \label{v2}
\end{equation}
Notice that the definition of $f\otimes_r g$ does not depend
on the particular choice of $\{e_k,\,k\geq 1\}$, and that
$f\otimes _{r}g$ is not necessarily symmetric; we denote its
symmetrization by $f\widetilde{\otimes }_{r}g\in \EuFrak H^{\odot (p+q-2r)}$.
Moreover, $f\otimes _{0}g=f\otimes g$ equals the tensor product of $f$ and
$g$ while, for $p=q$, $f\otimes _{q}g=\langle f,g\rangle _{\EuFrak H^{\otimes q}}$. 
When $\HH = L^2(A,\mathcal{A},\mu)$ and $r=1,...,p\wedge q$, the contraction $f\otimes _{r}g$ is the element of $L^2(\mu^{p+q-2r})$ given by
\begin{eqnarray}\label{e:contraction}
&& f\otimes _{r}g (x_1,...,x_{p+q-2r})\\
&& = \int_{A^r} f(x_1,...,x_{p-r},a_1,...,a_r)g(x_{p-r+1},...,x_{p+q-2r},a_1,...,a_r)d\mu(a_1)...d\mu(a_r). \notag
\end{eqnarray}

It can also be shown that the following {\sl multiplication formula} holds: if $f\in \EuFrak
H^{\odot p}$ and $g\in \EuFrak
H^{\odot q}$, then
\begin{eqnarray}\label{multiplication}
I_p(f) I_q(g) = \sum_{r=0}^{p \wedge q} r! {p \choose r}{ q \choose r} I_{p+q-2r} (f\widetilde{\otimes}_{r}g).
\end{eqnarray}
\smallskip

\subsection{Malliavin operators}\label{ss:mall}

We now introduce some basic elements of the Malliavin calculus with respect
to the isonormal Gaussian process $X$. Let $\mathcal{S}$
be the set of all
cylindrical random variables of
the form
\begin{equation}
F=g\left( X(\phi _{1}),\ldots ,X(\phi _{n})\right) ,  \label{v3}
\end{equation}
where $n\geq 1$, $g:\mathbb{R}^{n}\rightarrow \mathbb{R}$ is an infinitely
differentiable function such that its partial derivatives have polynomial growth, and $\phi _{i}\in \EuFrak H$,
$i=1,\ldots,n$.
The {\it Malliavin derivative}  of $F$ with respect to $X$ is the element of
$L^2(\Omega ,\EuFrak H)$ defined as
\begin{equation*}
DF\;=\;\sum_{i=1}^{n}\frac{\partial g}{\partial x_{i}}\left( X(\phi
_{1}),\ldots ,X(\phi _{n})\right) \phi _{i}.
\end{equation*}
In particular, $DX(h)=h$ for every $h\in \EuFrak H$. By iteration, one can
define the $m$th derivative $D^{m}F$, which is an element of $L^2(\Omega ,\EuFrak H^{\odot m})$,
for every $m\geq 2$.
For $m\geq 1$ and $p\geq 1$, ${\mathbb{D}}^{m,p}$ denotes the closure of
$\mathcal{S}$ with respect to the norm $\Vert \cdot \Vert _{m,p}$, defined by
the relation
\begin{equation*}
\Vert F\Vert _{m,p}^{p}\;=\;E\left[ |F|^{p}\right] +\sum_{i=1}^{m}E\left[
\Vert D^{i}F\Vert _{\EuFrak H^{\otimes i}}^{p}\right].
\end{equation*}
We often use the notation $\mathbb{D}^{\infty} := \bigcap_{m\geq 1}
\bigcap_{p\geq 1}\mathbb{D}^{m,p}$.

\begin{rem}\label{r:density}{\rm Any random variable $Y$ that is a finite linear combination of multiple Wiener-It\^o integrals is an element of $\mathbb{D}^\infty$. 
Moreover, if $Y\neq 0$, then the law of $Y$ admits a density with respect to the Lebesgue measure -- see \cite{Shigekawa}.
}
\end{rem}

The Malliavin derivative $D$ obeys the following \textsl{chain rule}. If
$\varphi :\mathbb{R}^{n}\rightarrow \mathbb{R}$ is continuously
differentiable with bounded partial derivatives and if $F=(F_{1},\ldots
,F_{n})$ is a vector of elements of ${\mathbb{D}}^{1,2}$, then $\varphi
(F)\in {\mathbb{D}}^{1,2}$ and
\begin{equation}\label{e:chainrule}
D\,\varphi (F)=\sum_{i=1}^{n}\frac{\partial \varphi }{\partial x_{i}}(F)DF_{i}.
\end{equation}

\begin{rem}\label{r:polg}{\rm
By approximation, it is easily checked that equation (\ref{e:chainrule}) continues to hold in the following two cases:
(i) $F_i \in \mathbb{D}^\infty$ and $\varphi$ has continuous partial derivatives with at most polynomial growth, and
(ii) $F_i\in \mathbb{D}^{1,2}$ has an absolutely continuous distribution and $\varphi$ is Lipschitz continuous.
}
\end{rem}

Note also that a random variable $F$ as in (\ref{E}) is in ${\mathbb{D}}^{1,2}$ if and only if
$\sum_{q=1}^{\infty }q\|J_qF\|^2_{L^2(\Omega)}<\infty$
and, in this case, $E\left[ \Vert DF\Vert _{\EuFrak H}^{2}\right]
=\sum_{q=1}^{\infty }q\|J_qF\|^2_{L^2(\Omega)}$. If $\EuFrak H=
L^{2}(A,\mathcal{A},\mu )$ (with $\mu $ non-atomic), then the
derivative of a random variable $F$ as in (\ref{E}) can be identified with
the element of $L^2(A \times \Omega )$ given by
\begin{equation}
D_{x}F=\sum_{q=1}^{\infty }qI_{q-1}\left( f_{q}(\cdot ,x)\right) ,\quad x \in A.  \label{dtf}
\end{equation}

We denote by $\delta $ the adjoint of the operator $D$, also called the
\textit{divergence operator}. A random element $u\in L^2(\Omega ,\EuFrak H)$
belongs to the domain of $\delta $, noted $\mathrm{Dom}\,\delta $, if and
only if it verifies
$|E\langle DF,u\rangle _{\EuFrak H}|\leq c_{u}\,\Vert F\Vert _{L^2(\Omega)}$
for any $F\in \mathbb{D}^{1,2}$, where $c_{u}$ is a constant depending only
on $u$. If $u\in \mathrm{Dom}\,\delta $, then the random variable $\delta (u)$
is defined by the duality relationship (customarily called \textit{integration by parts
formula})
\begin{equation}
E[F\delta (u)]=E[\langle DF,u\rangle _{\EuFrak H}],  \label{ipp0}
\end{equation}
which holds for every $F\in {\mathbb{D}}^{1,2}$.

The operator $L$, defined as
$L=\sum_{q=0}^{\infty }-qJ_{q}$,
is the {\it infinitesimal generator of the Ornstein-Uhlenbeck
semigroup}.
The domain of $L$ is
\begin{equation*}
\mathrm{Dom}L=\{F\in L^2(\Omega ):\sum_{q=1}^{\infty }q^{2}\left\|
J_{q}F\right\| _{L^2(\Omega )}^{2}<\infty \}=\mathbb{D}^{2,2}\text{.}
\end{equation*}
There is an important relation between the operators $D$, $\delta $ and $L$.
A random variable $F$ belongs to
$\mathbb{D}^{2,2}$ if and only if $F\in \mathrm{Dom}\left( \delta D\right) $
(i.e. $F\in {\mathbb{D}}^{1,2}$ and $DF\in \mathrm{Dom}\delta $) and, in
this case,
\begin{equation}
\delta DF=-LF.  \label{k1}
\end{equation}

For any $F \in L^2(\Omega )$, we define $L^{-1}F =\sum_{q=1}^{\infty }-\frac{1}{q} J_{q}(F)$. The operator $L^{-1}$ is called the
\textit{pseudo-inverse} of $L$. Indeed, for any $F \in L^2(\Omega )$, we have that $L^{-1} F \in  \mathrm{Dom}L
= \mathbb{D}^{2,2}$,
and
\begin{equation}\label{Lmoins1}
LL^{-1} F = F - E(F).
\end{equation}

The following result is used throughout the paper.

\begin{lemme}\label{L : Tech1}
Suppose that $H\in\mathbb{D}^{1,2}$ and $G\in L^2(\Omega)$. Then,
$L^{-1}G \in \mathbb{D}^{2,2}$ and
\begin{equation}
E[HG] = E[H]E[G]+E[\langle DH,-DL^{-1}G\rangle_{\HH}].
\end{equation}
\end{lemme}
{\it Proof}.
By (\ref{k1}) and (\ref{Lmoins1}),
\[
E[HG]-E[H]E[G]=E[H(G-E[G])]=E[H\times LL^{-1}G]=E[H \delta(- DL^{-1}G) ],
\]
and the result is obtained by using the integration by parts formula (\ref{ipp0}).
\fin

\section{Stein's equations and cumulants}\label{s:steincum}
In order to prove our main results, we shall combine the integration by parts formula of Malliavin calculus, both with a standard version of the Stein's method for normal approximations 
(see \cite{ChenGoldShao} for an exhaustive presentation of this technique) and with a fine analysis of the cumulants associated with random variables living in a fixed chaos. One of our main tools is an Edgeworth-type expansion (inspired by Barbour's paper \cite{BarbourPtrf1986}) for smooth transformations of Malliavin differentiable random variables. These fundamental topics are presented in the three subsections to follow.

\subsection{Stein's equations and associated bounds}\label{ss:stein}

Let $N \sim \mathscr{N}(0,1)$ be a standard Gaussian random variable, and let 
$h : \R \rightarrow \R$ be a continuous 
function. 
\begin{rem}
{\rm
In the literature about Stein's method and normal approximation, it is customary at this stage to assume
that $h$ is merely Borel measurable. However, this leads to some technical issues that are not necessary here. See e.g. \cite[Chapter 3]{np-book}.
}
\end{rem}

We associate with $h$ the following {\it Stein's equation}:
\begin{equation}\label{e:steinEQ}
h(x) - E[h(N)] = f'(x) - xf(x), \quad x\in \R.
\end{equation}
It is easily checked that, if $E|h(N)|<\infty$, then the function
\begin{equation}\label{fh}
f_h(x) = e^{x^2/2} \int_{-\infty}^x (h(y) - E[h(N)])e^{-y^2/2} dy, \quad x\in\R,
\end{equation}
is the unique solution of (\ref{e:steinEQ}) verifying the additional asymptotic condition \[\lim_{x\rightarrow \pm \infty} f_h(x)e^{-x^2/2} = 0.\]

\smallskip

In this paper, we will actually deal with Stein's equations associated with functions $h$ that are differentiable up to a certain order. The
following statement (proved by Daly in \cite{daly}) is an important tool for our analysis. Throughout the following, given a
smooth function $g: \R\rightarrow\R$, we shall denote by $g^{(k)}$, $k=1,2,...$, the $k$th derivative of $g$; we sometimes write $g' = g^{(1)}$, $g'' = g^{(2)}$, and so on.

\begin{prop}\label{p:daly} Let the previous notation prevail, fix an integer $k\geq 0$, and assume that the function
$h$ is $(k+1)$-times differentiable and such that $h^{(k)}$ is absolutely continuous. Then, $f_h$ is $(k+2)$-times differentiable,
and one has the estimate
\begin{equation}\label{e:daly}
\|f_h^{(k+2)}\|_\infty \leq 2\|h^{(k+1)}\|_\infty.
\end{equation}
Moreover, the continuity of $h^{(k+1)}$ implies the continuity of $f_h^{(k+2)}$.
\end{prop}
{\it Proof}. The first part, i.e., inequality (\ref{e:daly}), is exactly Theorem 1.1 of \cite{daly}, whereas
the transfer of continuity is easily checked
by induction and by using (\ref{e:steinEQ}).\fin

\subsection{Cumulants}\label{ss:cum}

We now formally define the cumulants associated with a random variable.

\begin{defi}[Cumulants]\label{D : cum}{\rm Let $F$ be a real-valued random variable such that $E|F|^m<\infty$ for some integer
$m\geq 1$, and define $\phi_F(t) = E[e^{itF}]$, $t\in\R$, to be the characteristic function of $F$.
Then, for $j=1,...,m$, the $j$th {\it cumulant} of $F$, denoted by $\kappa_j(F)$, is given by
\begin{equation}
\kappa_j (F) = (-i)^j \frac{d^j}{d t^j} \log \phi_F (t)|_{t=0}.
\end{equation}}
\end{defi}

\begin{rem}{\rm The first four cumulants are the following: $\kappa_1(F) = E[F]$, $\kappa_2(F) = E[F^2]-E[F]^2 = {\rm Var}(F)$,
$\kappa_3(F) = E[F^3]-3E[F^2]E[F]+2E[F]^3$, and
\[
\kappa_4(F) = E[F^4] - 3E[F]E[F^3]-3E[F^2]^2+12E[F]^2E[F^2]-6E[F]^4.
\]
In particular, when $E[F]=0$ one sees that $\kappa_3(F) = E[F^3]$ and $\kappa_4(F) = E[F^4]-3E[F^2]^2$.}
\end{rem}

\smallskip

The reader is referred to \cite[Chapter 3]{PeTa} for a self-contained presentation of the properties of cumulants and for several
combinatorial characterizations. The following relation (proved e.g. in \cite[Proposition 2.2]{np-jfa}) shows that moments can be
recursively defined in terms of cumulants (and vice-versa): fix $m= 1,2...$, and assume that $E|F|^{m+1}<\infty$, then
\begin{equation}\label{EQ : RecMom}
E[F^{m+1}] = \sum_{s=0}^m \binom{m}{s}\kappa_{s+1}(F) E[F^{m-s}].
\end{equation}

We now want to characterize cumulants in terms of Malliavin operators. To do so, we need the following recursive definition (taken from \cite{np-jfa}).

\begin{defi}\label{Def : Gamma}{\rm Let $F\in \mathbb{D}^{\infty}$. The sequence of random variables $\{\Gamma_j(F) : j\geq 0\}\subset
\mathbb{D}^\infty$ is recursively defined as follows. Set $\Gamma_0(F) = F$
and, for every $j\geq 1$, \[
\Gamma_{j}(F) = \langle DF,-DL^{-1}\Gamma_{j-1}(F)\rangle_{\HH}.
\]
Note that each $\Gamma_j(F)$ is a well-defined element of $ \mathbb{D}^\infty$, since $F$ is assumed to be in $\mathbb{D}^\infty$ -- see \cite[Lemma 4.2(3)]{np-jfa}}
\end{defi}

For instance, one has that $\Gamma_1(F) = \langle DF,-DL^{-1}F\rangle_{\HH}$. The following statement provides two explicit relations
((\ref{e:cumgamma}) and (\ref{e:cumgamma2})) connecting the random variables $\Gamma_j(F)$ to the cumulants of $F$. Equation (\ref{e:cumgamma}) has been proved in \cite[Theorem 4.3]{np-jfa}, whereas (\ref{e:cumgamma2}) is new.

\begin{prop} Let $F\in \mathbb{D}^\infty$. Then $F$ has finite moments of every order, and the following relation holds for every $s\geq 0$:
\begin{equation}\label{e:cumgamma}
\kappa_{s+1}(F) = s!E[\Gamma_s(F)].
\end{equation}
If moreover $E(F)=0$ then, for every $s\geq 1$,
\begin{equation}\label{e:cumgamma2}
\kappa_{s+2} (F)= \frac{1}{2}(s+1)!\,E\left[F^2\left(\Gamma_{s-1}(F) - \frac{\kappa_s(F)}{(s-1)!}\right)\right].
\end{equation}
\end{prop}
{\it Proof.} In view of \cite[Theorem 4.3]{np-jfa}, we have only to prove (\ref{e:cumgamma2}). Applying Lemma \ref{L : Tech1} in the special case $H = F^2$ and $G= \Gamma_{s-1}(F)$, and using the relation $DF^2 = 2FDF$, one deduces that
\[
E[F^2\Gamma_{s-1}(F)] = E[F^2]E[\Gamma_{s-1}(F)] +2E[F \Gamma_s(F)].\]
Now apply Lemma \ref{L : Tech1} in the case $H=F$ and $G=\Gamma_s(F)$: exploiting the fact that $F$ is centered together with (\ref{e:cumgamma}), we infer that
\begin{equation}\label{azertiop}
E[F \Gamma_s(F)] =E[F( \Gamma_s(F)-E[\Gamma_s(F)])] =  E[ \Gamma_{s+1}(F)] =\frac{\kappa_{s+2}(F)}{(s+1)!}.
\end{equation}
Since (\ref{e:cumgamma}) implies that $(s-1)!E[\Gamma_{s-1}(F)]  = \kappa_s(F)$, the conclusion follows.
\fin

\begin{rem}{\rm

\begin{enumerate}

\item Relation (\ref{e:cumgamma}) continues to hold under weaker assumptions on the regularity of $F$. See again \cite[Theorem 4.3]{np-jfa}.

\item Relation (\ref{e:cumgamma}) generalizes the following well-known fact: if $F\in\mathbb{D}^{1,2}$, then
$\Gamma_1(F) = \langle DF,-DL^{-1}F\rangle_\HH$ is in $L^1(\Omega)$ and
\begin{equation}\label{e:var}
{\rm Var}(F) = E[\Gamma_1(F)].
\end{equation}

\end{enumerate}
}
\end{rem}

The following statement provides an explicit expression for $\Gamma_s(F)$, $s\geq 1$, when
$F$ has the form of a multiple integral.

\begin{prop}[see \cite{np-jfa}, formula (5.25)]
\label{thm-pasmal}
Let $q\geq 2$, and assume that $F=I_q(f)$ with $f\in\HH^{\odot q}$.
Then, for any $s\geq 1$, we have
\begin{eqnarray}
\Gamma_{s}(F)&=&
\sum_{r_1=1}^{q} \ldots\sum_{r_{s}=1}^{[sq-2r_1-\ldots-2r_{s-1}]\wedge q}c_q(r_1,\ldots,r_{s})
{\bf 1}_{\{r_1< q\}}
\ldots {\bf 1}_{\{
r_1+\ldots+r_{s-1}< \frac{sq}2
\}}\label{for2}\\
&&
\hskip5cm\times I_{(s+1)q-2r_1-\ldots-2r_{s}}\big(
(...(f\widetilde{\otimes}_{r_1} f) \widetilde{\otimes}_{r_2} f)\ldots
f)\widetilde{\otimes}_{r_{s}}f
\big),\notag
\end{eqnarray}
where the constants $c_q(r_1,\ldots,r_{s-2})$ are recursively defined as follows:
\[
c_q(r)=q(r-1)!\binom{q-1}{r-1}^2,
\]
and, for $a\geq 2$,
\[
c_q(r_1,\ldots,r_{a})=q(r_{a}-1)!\binom{aq-2r_1-\ldots - 2r_{a-1}-1}{r_{a}-1}
\binom{q-1}{r_{a}-1}c_q(r_1,\ldots,r_{a-1}).
\]
\end{prop}

\smallskip

\begin{rem}{\rm
By combining (\ref{e:cumgamma}) with (\ref{for2}), we immediately get a representation of cumulants that is alternative to the one based on
`diagram formulae'. See \cite[Theorem 5.1]{np-jfa} for details on this point.
}
\end{rem}

\subsection{Assessing Edgeworth-type expansions}\label{ss:edge}
The following Edgeworth-type expansion also plays a crucial role in the following.

\begin{prop}\label{p:notrebarbour}
Let $F$ be an element of $\mathbb{D}^{\infty}$. Then, for every $M\geq 1$ and every function $f : \mathbb{R}\rightarrow \R$ that is
$M$ times continuously differentiable with derivatives having at most polynomial growth, we have
\begin{equation}\label{e:notrebarbour}
E[Ff(F)] = \sum_{s=0}^{M-1} \frac{\kappa_{s+1}(F)}{s!} E[f^{(s)}(F)] + E[\Gamma_{M}(F)f^{(M)}(F)].
\end{equation}
\end{prop}
{\it Proof.} Using twice Lemma \ref{L : Tech1}, first in the case $H = F$ and $G = f(F)$ and then in the case $F = \Gamma_1(F)$ and $G = f'(F)$, we 
deduce that
\begin{eqnarray*}
E[Ff(F)] &=& E[F]E[f(F)]+E[f'(F)\Gamma_1(F)] \\
&=& E[F]E[f(F)]+E[f'(F)]E[\Gamma_1(F)] +E[f''(F)\Gamma_2(F)],
\end{eqnarray*}
where we have used the chain rule (\ref{e:chainrule}) as well as Remark \ref{r:polg}. 
Therefore, (\ref{e:notrebarbour}) holds for $M=1,2$ (see also (\ref{e:var})).
The case of a general $M$ follows immediately from an induction argument and by using (\ref{e:cumgamma}).
\fin

The following statements contain two important consequences of (\ref{e:notrebarbour}). They will be used in order to prove our main findings.

\begin{cor}\label{c:est1}
Let $N\sim \mathscr{N}(0,1)$ and fix $F\in \mathbb{D}^\infty$ such that $E[F]=0$, $E[F^2] = \kappa_2(F) = 1$. For $M\geq 2$ , let $h : \R\rightarrow\R$ be $(M-1)$ times continuously differentiable with bounded derivatives, and define $f_h$ according to (\ref{fh}). Then,
\begin{equation}\label{PBB1}
\left| E[h(N)] - E[h(F)] - \sum_{s=2}^{M-1} \frac{\kappa_{s+1}(F)}{s!} E[f_h^{(s)}(F)]  \right| \leq 2\|h^{(M-1)}\|_\infty E|\Gamma_{M}(F)|.
\end{equation}
\end{cor}
{\it Proof.} From Proposition \ref{p:daly}, we deduce that the function $f_h$ is $M$-times continuously differentiable
and that, for $k=2,...,M$,  $ \|f_h^{(k)}\|_\infty\leq 2\|h^{(k-1)}\|_\infty$.
Using a Taylor expansion, we deduce that $f'_h$ has at most polynomial growth. It follows that (\ref{e:notrebarbour})
can be applied to the function $f_h$, and the conclusion is obtained from the relation $E[h(N)] - E[h(F)] = E[Ff_h(F)] - E[f'_h(F)]$.
\fin

\begin{cor}\label{c:est2}
Let $N\sim \mathscr{N}(0,1)$ and fix $F\in \mathbb{D}^\infty$ such that $E[F]=0$ and $E[F^2] = \kappa_2(F) = 1$. Let  $h : \R\rightarrow\R$ be twice continuously differentiable and such that $\|h''\|_\infty \leq 1$, and define $f_h$ according to (\ref{fh}). Then,
\begin{equation}\label{PBB2}
\big| E[h(F)]-E[h(N)]\big|\leq K |E[F^3]|+2 E|\Gamma_3(F)|;
\end{equation}
where $K: = 1+ E[|F|]$.
\end{cor}
{\it Proof.} We first observe that
$E[h(F)]-E[h(N)]=E[\tilde{h}(F)]-E[\tilde{h}(N)]$, where
$\tilde{h}(x)=h(x)-h(0)-h'(0)x$, so that we can assume without loss of
generality that $h(0)=h'(0)=0$. Thus, because $\|h''\|_\infty\leq 1$, we
get that $|h(x)|\leq \frac{x^2}{2}$ and $|h'(x)|\leq |x|$ for all
$x\in\R$, while $|E[h(N)]|\leq \frac 12$. It
follows from \eqref{fh} that
\[
|f_h(0)|\leq\int_{0}^\infty
\left(\frac{y^2}{2}+\frac12\right)e^{-y^2/2}dy=
\frac{\sqrt{2\pi}}2\leq
2.
\]
Next, Proposition \ref{p:daly} shows that $f_h$ is thrice continuously
differentiable with
$\|f_h'''\|_\infty\leq 2\,\|h''\|_\infty\leq 2$.
On the other hand, for all $x\in\R$,
\begin{eqnarray*}
f'_h(x)&=&xf_h(x) + h(x) - E[h(N)],\\
f''_h(x)&=&f_h(x) + x f'_h(x) + h'(x).
\end{eqnarray*}
Consequently, $f''_h(0)=f_h(0)$ and
\begin{equation}|f''_h(x)|\leq |f_h(0)|+|f''_h(x)-f''_h(0)|\leq 2+\|f_h'''\|_\infty|x|\leq 2+2|x|.\label{f_h}
\end{equation}
We deduce that $|E[f''_h(F)]|\leq 2K$.
Applying (\ref{e:notrebarbour}) to $f_h$ in the case $M=3$ yields therefore
\begin{equation}\label{expa}
E[Ff_h(F)]=E[f'_h(F)] + \frac12E[f''_h(F)]E[F^3]+E[f_h'''(F)\Gamma_3(F)],
\end{equation}
implying in turn that
\[
\big| E[h(F)]-E[h(N)]\big|\leq \frac12|E[f''_h(F)]| |E[F^3]|+|f_h'''|_\infty E|\Gamma_3(F)|,
\]
from which the desired conclusion follows.

\fin

\begin{rem}{\rm
\begin{enumerate}
\item The idea of bounding quantities of the type
\[
\left| E[Ff(F)]  - \sum_{s=0}^{M-1} \frac{\kappa_{s+1}(F)}{s!} E[f^{(s)}(F)]\right|,
\]
in order to estimate the distance between $F$ and $N\sim\mathscr{N}(0,1)$, dates back to Barbour's seminal paper
\cite{BarbourPtrf1986}. Due to the fact that $F$ is a smooth functional of a Gaussian field,
observe that the expression of the `rest'
$E[\Gamma_{M}(F)f^{(M)}(F)] $ appearing in (\ref{e:notrebarbour}) is remarkably simpler than the ones
computed in \cite{BarbourPtrf1986}.

\item For a fixed $M$, the expansion (\ref{e:notrebarbour}) may hold under weaker assumptions on $F$ and $f$. For instance,
if $F\in \mathbb{D}^{1,2}$ has an absolutely continuous law, then, for every Lipschitz continuous function $f:\R\rightarrow \R$,
\begin{equation}\label{q}
E[Ff(F)] = E[F]E[f'(F)]+E[f'(F)\Gamma_1(F)].
\end{equation}
Equation (\ref{q}) is the starting point of the analysis developed in \cite{np-ptrf}.
\end{enumerate}
}
\end{rem}

\section{Some technical estimates}\label{s:estimates}
This section contains several estimates that are needed in the proof of Theorem \ref{t:mainupper} and Theorem \ref{t:mainlower}.

\subsection{Inequalities for kernels}

For every integer $M\geq 1$, we write $[M] = \{1,...,M\}$. Fix a set $Z$, as well as a vector ${\bf z} = (z_1,...,z_M)\in Z^M$ and a nonempty set $b\subseteq [M]$: we denote by ${\bf z}_b$ the element of $Z^{|b|}$ (where $|b|$ is the cardinality of $b$) obtained by deleting from ${\bf z}$ the entries whose index is not contained in $b$. For instance, if $M=5$ and $b=\{1,3,5\}$, then ${\bf z}_b = (z_1,z_3,z_5)$. Now consider the following setting:

\begin{enumerate}

\item[{\bf ($\alpha$)}] $(Z,\mathcal{Z})$ is a measurable space, and $\mu$ is a measure on it;

\item[{\bf ($\beta$)}] $B,q\geq 2$ are integers, and $b_1,...,b_q$  are nonempty subsets of $[B]$ such that $\cup_{i} b_i = [B]$, and each $ k \in [B]$ appears in exactly two of the $b_i$`s (this implies in particular that $\sum_i |b_i| = 2B$, and also that, if $q=2$, then necessarily $b_1=b_2 = [B]$);

\item[{\bf ($\gamma$)}] $F_1,...,F_q$ are functions such that $F_i \in L^2(Z^{|b_i|},\mathcal{Z}^{|b_i|},\mu^{|b_i|})=L^2(\mu^{|b_i|})$ for every $i=1,...,q$ (in particular, each $F_i$ is a function of $|b_i|$ variables).

\end{enumerate}

The following generalization of the Cauchy-Schwarz inequality is crucial in this paper.

\begin{lemme}[Generalized Cauchy-Schwarz Inequality]\label{l:CS}

Under assumptions {\bf ($\alpha$)}-{\bf ($\beta$)}-{\bf ($\gamma$)}, the following inequality holds:
\begin{equation}\label{e:genCS}
\int_{Z^B} \prod_{i=1}^q |F_i({\bf z}_{b_i})| \mu(dz_1)\cdot\cdot\cdot\mu(dz_B) \leq \prod_{i=1}^q \|F_i\|_{L^2(\mu^{|b_i|})}.
\end{equation}
\end{lemme}
{\it Proof}.
The case $q=2$ is just the Cauchy-Schwarz inequality, and the general result is obtained by recursion on $q$. The argument goes as follows: call $A$ the left-hand side of (\ref{e:genCS}), and assume that the
desired estimate is true for $q-1$. Applying the Cauchy-Schwarz inequality, we deduce that (with obvious notation)
\[
A\leq \|F_1\|_{L^2(\mu^{|b_1|})} \times \left( \int_{Z^{|b_1|}} \Phi({\bf z}_{b_1})^2 \mu^{|b_1|}(d{\bf z}_{b_1})\right)^{1/2},
\]
where the quantity $\Phi({\bf z}_{b_1})$ is obtained by integrating the product $\prod_{i=2}^q |F_i({\bf z}_{b_i})|$  over those variables $z_j$ such that $j \not\in b_1$. More explicitly, writing $J$ for the class of those $i \in \{2,...,q\}$ such that $b_i \subseteq b_1$,
\begin{equation}\label{z}
\Phi({\bf z}_{b_1}) =\prod_{i\in J} |F_i({\bf z}_{b_i})| \times \int_{Z^{|B| - |b_1|}} \prod_{i\in J^c} |F_i({\bf z}_{b_i\cap b_1},\, {\bf z}_{b_i\cap b^c_1} )| \mu^{|B| - |b_1|}(d{\bf z}_{b_1^c}),
\end{equation}
where $b_1^c$ and $J^c$ indicate, respectively, the complement of $b_1$ (in $[B]$) and the complement of $J$ (in $ \{2,...,q\}$), and $\int \prod_{j\in \emptyset} = 1$ by convention.  By construction, one has that the sets $b_i$ such that $i\in J$ are disjoint, and also that $b_i\cap b_j = \emptyset$, for every $i\in J$ and $j\in J^c$. If $J^c = \emptyset$, there is nothing to prove. If $J^c\neq \emptyset$, one has to observe that the blocks $b'_i = b_i \cap b_1^c$, $i\in J^c$, verify assumption {\bf ($\beta$)} with respect to the set $ [B]\backslash b_1$ (that is, the class $\{b'_i : i\in J^c\}$ is composed of nonempty subsets of $[B]\backslash b_1$ such that $\cup_{i} b'_i =[B]\backslash b_1$, and each $ k \in [B]\backslash b_1$ appears in exactly two of the $b'_i$'s). Since $|J^c| \leq q-1$, the recurrence assumption can be applied to the integral on the right-hand side of (\ref{z}), thus yielding the desired conclusion.
\fin

\medskip

Let $\HH$ be a real separable Hilbert space. The next estimates also play a pivotal role in our arguments.

\begin{lemme}
Let $p,q\geq 1$ be two integers, $r\in\{0,\ldots,p\wedge q\}$, and $u\in\HH^{\odot p}$, $v\in\HH^{\odot q}$
Then
\begin{eqnarray}
\|u\,\,\widetilde{\otimes}_r \,v\|_{\HH^{\otimes(p+q-2r)}}&\leq& \|u\otimes_r v\|_{\HH^{\otimes(p+q-2r)}}\label{1}\\
\|u\otimes_r v\|
_{\HH^{\otimes(p+q-2r)}
}&\leq&  \|u\|_{\HH^{\otimes p}}\sqrt{\|v\otimes_{q-r}v\|_{\HH^{\otimes (2r)}}}\leq \|u\|_{\HH^{\otimes p}}\|v\|_{\HH^{\otimes q}}\label{2}.
\end{eqnarray}
Moreover, if $q!\|v\|^2_{\HH^{\otimes q}}=1$ (that is, if $E[I_q(v)^2]=1$), then
\begin{equation}\label{3}
\max_{1\leq r\leq q-1}\|v\otimes_r v\|^2_{\HH^{\otimes 2q-2r}}\leq \frac{\kappa_4(I_q(v))}{q!^2q^2}.
\end{equation}
\end{lemme}
{\it Proof}. The proof of (\ref{1}) is evident, by using the very definition of a symmetrized function.
To show the first inequality in (\ref{2}), apply first Fubini to get that
$
\|u\otimes_r v\|^2
_{\HH^{\otimes(p+q-2r)}} = \langle u\otimes_{p-r} u, v\otimes_{q-r}v\rangle_{\HH^{\otimes(2r)}},
$
and then Cauchy-Schwarz to get the desired conclusion. The second inequality in (\ref{2}) is an immediate consequence
of Cauchy-Schwarz. Finally, the proof of (\ref{3}) is obtained by using \cite[first equality on p. 183]{nunugio}.
\fin

\subsection{Inequalities for cumulants and related quantities}

The following proposition contains all the estimates that are necessary for proving the main results in the paper.
For every random variable $Y$ such that $E|Y|^m<\infty$ ($m\geq 1$), we denote by $\kappa_m(Y)$ the $m$th cumulant of $Y$ -- see
Definition \ref{D : cum}. Given a vector ${\bf z} = (z_1,...,z_d)$ and a permutation $\sigma$ of $[d]$, we write $\sigma({\bf z}) =
(z_{\sigma(1)}, ...,z_{\sigma(d)})$. Given a function $F(z_1,...,z_d)$ of $d$ variables and a permutation $\sigma$ of $[d]$, we write
\[
(F)_\sigma({\bf z}) =F(\sigma({\bf z}))= F(z_{\sigma(1)},...,z_{\sigma(d)}).
\]
Also, for vectors ${\bf z} = (z_1,...,z_j)$ and ${\bf y} = (y_1,...,y_k)$, we shall write ${\bf z}\vee {\bf y}$ for the vector of dimension 
$j+k$ obtained by juxtaposing  ${\bf z}$ and ${\bf y}$, that is,  ${\bf z}\vee {\bf y} = (z_1,...,z_j,y_1,...,y_k)$. Finally, in the subsequent 
proofs we will identify vectors of dimension zero with the empty set: if ${\bf z}$ has dimension zero, then integration with respect to 
${\bf z}$ is removed by convention.
\begin{prop}\label{prop-gamma} 
We use the notation introduced in Definitions \ref{D : cum} and \ref{Def : Gamma}.
For each integer $q\geq 2$ there exists positive constants $c_2(q),c_3(q),c_4(q) $ (only depending on $q$)
such that, for all $F=I_q(f)$ with $f\in\HH^{\odot q}$ and $E[F^2]=1$,
we have
\begin{eqnarray}
E\left[\left|\Gamma_2(F)-\frac12\kappa_3(F)\right|\right]& \leq & c_2(q)\times\kappa_4(F)^{\frac34} \label{E0},\\
E[|\Gamma_3(F)|]& \leq & c_3(q)\times\kappa_4(F) \label{E1},\\
E[|\Gamma_4 (F)|] &\leq &     c_4(q) \times\kappa_4(F)^{\frac54}\label{E2}.
\end{eqnarray}
\end{prop}
\noindent {\it Proof}. By (\ref{e:cumgamma}), we have $s!E(\Gamma_s) =\kappa_{s+1}(F)$ for every $s\geq 1$.
Moreover, when $F = I_q(f)$ is as in the statement, recall the following explicit representation (\ref{for2}):
\begin{eqnarray}
\Gamma_{s}(F)&=&
\sum_{r_1=1}^{q} \ldots\sum_{r_{s}=1}^{[sq-2r_1-\ldots-2r_{s-1}]\wedge q}c_q(r_1,\ldots,r_{s})
{\bf 1}_{\{r_1< q\}}
\ldots {\bf 1}_{\{
r_1+\ldots+r_{s-1}< \frac{sq}2
\}}\label{for}\\
&&
\hskip5cm\times I_{(s+1)q-2r_1-\ldots-2r_{s}}\big(
(...(f\widetilde{\otimes}_{r_1} f) \widetilde{\otimes}_{r_2} f)\ldots
f)\widetilde{\otimes}_{r_{s}}f
\big).\notag
\end{eqnarray}
Without loss of generality, throughout the proof we shall assume that $\HH = L^2(Z,\mathcal{Z},\mu)$, where $Z$ is a Polish space, $\mathcal{Z}$ is the associated Borel $\sigma$-field, and $\mu$ is a $\sigma$-finite measure.

\medskip
\noindent{\it Proof of (\ref{E0})}.
According to (\ref{e:cumgamma}), one has that $E[\Gamma_2(F)] = \frac12\kappa_3(F)$, so that the random variable $E[\Gamma_2(F)] - \frac12\kappa_3(F)$ is obtained by restricting the sum in (\ref{for}) (in the case $s=2$) to the terms 
such that $r_1+r_2<\frac{3q}2$. By virtue of (\ref{3}), the inequality (\ref{E0}) will follow once it is shown that, for any choice of integers $r_1, r_2$ verifying such a constraint,

\begin{equation}\label{goal0}
\|
((f
\widetilde{\otimes}_{r_1}
f)
\widetilde{\otimes}_{r_2}
f)
\|_{\HH^{\otimes(3q-2r_1-2r_2)}}\leq \max_{1\leq r\leq q-1}\|f\otimes_r f\|^{\frac 32}_{\HH^{\otimes 2q-2r}}.
\end{equation}
Let us first assume that  $r_2<q$. Then $r_1$   and $q-r_2$ both belong to $\{1,\ldots,q-1\}$.
Thus, using  the two inequalities (\ref{1}) and  (\ref{2}), we infer that
\begin{eqnarray*}
&&\|
((f
\widetilde{\otimes}_{r_1}
f)
\widetilde{\otimes}_{r_2}
f)
\|_{\HH^{\otimes(3q-2r_1-2r_2)}}\\
&\leq &
\sqrt{\|f\otimes_{q-r_2}f\|_{\HH^{\otimes (2r_2)}}}
\,
\|
f
\otimes_{r_1}
f
\|_{\HH^{\otimes(2q-2r_1)}}\!
\leq\!
\max_{1\leq r\leq q-1}\!
\|
f
\otimes_{r}
f
\|_{\HH^{\otimes(2q-2r)}}^{\frac32}.
\end{eqnarray*}
Let us now consider the case when $r_2=q$ and $r_1<\frac q2$. The expression
\[
(f
\widetilde{\otimes}_{r_1}
f)
\widetilde{\otimes}_{q}
f
= \langle (f\widetilde{\otimes}_{r_1}f),f\rangle_{\HH^{\otimes q}}
\]
defines a function of $
q-2r_1$ variables. Taking into account the symmetry of $f$ and the symmetrization of contractions, such a function can be written as a convex linear combination of functions of the type
\begin{equation*}
 F({\bf t}) 
= \int f({\bf x}_1,{\bf t}_1,{\bf w})f({\bf x}_2,{\bf t}_2,{\bf w})f({\bf x}_1,{\bf x}_2)d\mu^{q+r_1}({\bf w}, {\bf x}_1,{\bf x}_2),
\end{equation*}
where ${\bf w}$ has length $r_1$, and ${\bf t}_1\vee {\bf t}_2 = \sigma({\bf t})$ for some permutation $\sigma$
and with ${\bf t}=(t_1,\ldots,t_{q-2r_1})$.
Without loss of generality we can assume that ${\bf t}_1$ has positive length
(recall that $r_1<q/2$ so that $q-2r_1>0$). We denote by $s_j$ the length of the vector ${\bf x}_j$. We 
then have  $1\leq s_1< q-r_1$ and $r_1< s_2\leq q-1$.  Exchanging the order of integrations, we can write
\begin{equation*}
 F({\bf t})
= \int f({\bf x}_1,{\bf t}_1,{\bf w})\left(f{\otimes}_{s_2}f\right)({\bf x}_1,{\bf t}_2,{\bf w})d\mu^{r_1+s_1}({\bf w}, {\bf x}_1).
\end{equation*}
Squaring $F$ and integrating, one sees that
\[
\|F\|^2_{\HH^{\otimes(q-2r_1)}} = \int \prod_{i=1}^3 (f\otimes_{\tau_i}f)_{\sigma_i}({\bf z}_{b_i})d\mu^B({\bf z}_{b_1},{\bf z}_{b_2},{\bf z}_{b_3}),
\]
with two of the $\tau_i$'s equal to $s_2$ and one to $q-r_1-s_1$, where $B=q+2s_1$, $\sigma_i$, $i=1,2,3,4$, is a permutation of $[2q-2\tau_i]$, and 
the sets $b_1,b_2,b_3$ verify property {\bf ($\beta$)}, as defined at the beginning of the present section. It follows from Lemma \ref{l:CS} that
\[
\|F\|_{\HH^{\otimes(q-2r_1)}} \leq  \max_{1\leq r\leq q-1}\|f\otimes_r f\|^{\frac 32}_{\HH^{\otimes 2q-2r}},
\]
from which we deduce (\ref{goal0}).

\smallskip

\noindent{\it Proof of (\ref{E1})}.
Our aim is to prove that for any  choice of $(r_1,r_2,r_3)$ appearing in the sum (\ref{for}) in the case $s=3$ one has the inequality
\begin{equation}\label{goal}
\|
((f
\widetilde{\otimes}_{r_1}
f)
\widetilde{\otimes}_{r_2}
f)
\widetilde{\otimes}_{r_3}
f
\|_{\HH^{\otimes(4q-2r_1-2r_2-2r_3)}}\leq \max_{1\leq r\leq q-1}\|f\otimes_r f\|^2_{\HH^{\otimes 2q-2r}}.
\end{equation}
 Remark that $((f\widetilde{\otimes}_{r_1}f)\widetilde{\otimes}_{r_2}f)$ has already been considered when looking at $\Gamma_2(F)-\frac12 \kappa_3(F)$, because of  the assumption that $r_1+r_2< \frac {3q}2$.

 So, using the previous estimates and \eqref{2}, we conclude directly for $r_3<q$. It remains to consider the case when $ r_3=q$.

As before, taking into account the symmetry of $f$ and the symmetrization of contractions, it is sufficient to consider functions of $2(q-r_1-r_2)$ variables of the type
\begin{eqnarray*}
&& F({\bf t}) = F(t_1,...,t_{2(q-r_1-r_2)})\\
&&= \int_{Z^{q+r_1+r_2}} \!\!\!\!\! f({\bf x}_1,{\bf a}_1,{\bf t}_1,{\bf w})f({\bf x}_2,{\bf a}_2,{\bf t}_2,{\bf w})\times \\
&&\quad\quad\quad\quad\quad\quad f({\bf a}_1,{\bf a}_2,{\bf t}_3,{\bf x}_3)f({\bf x}_1,{\bf x}_2,{\bf x}_3)\mu^{q+r_1+r_2}(d{\bf x}_1,d{\bf x}_2,d{\bf x}_3,d{\bf w},d{\bf a}_1,d{\bf a}_2),
\end{eqnarray*}
where ${\bf w}$ has length $r_1$, ${\bf a}_1\vee {\bf a}_2$ has length $r_2$ (with either ${\bf a}_1$ or ${\bf a}_2$ possibly equal to the empty set), and ${\bf t}_1\vee {\bf t}_2\vee{\bf t}_3 = \sigma({\bf t})$ for some permutation $\sigma$. Squaring $F$ and integrating, we claim that there exist integers $s_1,s_2,s_3,s_4 \in \{1,...,q-1\}$ such that
\[
\|F\|^2_{\HH^{2(q-r_1-r_2)}} = \int_{Z^{B}}\prod_{i=1}^4 (f\otimes_{s_i}f)_{\sigma_i}({\bf z}_{b_i})\mu(dz_1)\cdot\cdot\cdot\mu(dz_B),
\]
where $B=4q-2(s_1+ s_2+s_3+s_4)$, $\sigma_i$, $i=1,2,3,4$, is a permutation of $[2q-2s_i]$, and the sets $b_1,b_2,b_3,b_4$ verify property {\bf ($\beta$)}, as defined at the beginning of the present section. We have to consider separately two cases.

\smallskip

\noindent $(a)$: the length of ${\bf x}_3$ is not $0$: we can then take $s_1=s_2=r_1$ and $s_3=s_4$  equal to  the length of ${\bf x}_3$.

\noindent $(b)$: the length of ${\bf x}_3$ is $0$. Then either  ${\bf a}_1$ or ${\bf a}_2$ is not empty. Assuming that ${\bf a}_1$ is not empty, we can take for $s_1=s_2$ the length of ${\bf a}_1$ and for $s_3=s_4$ the length of ${\bf x}_2$, which is not $0$.

As before, it follows from Lemma \ref{l:CS} that
\[
\|F\|_{\HH^{2(q-r_1-r_2)}}\leq  \max_{1\leq r\leq q-1}\|f\otimes_r f\|^2_{\HH^{\otimes 2q-2r}},
\]
from which we deduce (\ref{goal}).

\smallskip

\noindent{\it Proof of (\ref{E2})}. Our aim is to prove that, for any  choice of $(r_1,r_2,r_3, r_4)$ which is present in the sum (\ref{for}) (in the case $s=4$),
we have
\begin{equation}\label{gogoal}
\|
(((f
\widetilde{\otimes}_{r_1}
f)
\widetilde{\otimes}_{r_2}
f)
\widetilde{\otimes}_{r_3}
f)\widetilde{\otimes}_{r_4} f
\|_{\HH^{\otimes(5q-2r_1-2r_2-2r_3-2r_4)}}\leq \max_{1\leq r\leq q-1}\|f\otimes_r f\|^{\frac 52}_{\HH^{\otimes 2q-2r}}.
\end{equation}

To do so, using the previous estimate \eqref{E1} and \eqref{2} we conclude directly for $r_4<q$. Hence, once again
it remains to consider the case when $ r_4=q$.

 As before, taking into account the symmetry of $f$ and the symmetrization of contractions, one has that the function $(((f
\widetilde{\otimes}_{r_1}
f)
\widetilde{\otimes}_{r_2}
f)
\widetilde{\otimes}_{r_3}
f)\widetilde{\otimes}_{r_4} f$ is a linear combination (with coefficients not depending on $f$) of functions in $3q-2r_1-2r_2-2r_3$ variables having the form
\begin{eqnarray*}
 F({\bf t}) &= &F(t_1,...,t_{3q-2r_1-2r_2-2r_3})\\
&=& \int_{Z^{q+r_1+r_2+r_3}} \!\!\!\!\! f({\bf x}_1,{\bf a}_1,{\bf b}_1,{\bf t}_1,{\bf w})f({\bf x}_2,{\bf a}_2,{\bf b}_2,{\bf t}_2,{\bf w})f({\bf x}_3,{\bf a}_1,{\bf a}_2,{\bf b}_3,{\bf t}_3)\times \\
&&\hfill\qquad\quad\quad\quad\quad\quad\quad f({\bf b}_1,{\bf b}_2,{\bf b}_3,{\bf t}_4,{\bf x}_4)f({\bf x}_1,{\bf x}_2,{\bf x}_3, {\bf x}_4)\mu^{q+r_1+r_2+r_3}(d{\bf x},d{\bf w},d{\bf a},d{\bf b}),
\end{eqnarray*}
where ${\bf w}$ has length $r_1$, ${\bf a}={\bf a}_1\vee {\bf a}_2$ has length $r_2$ (with either ${\bf a}_1$ or ${\bf a}_2$ possibly equal to the empty set), ${\bf b}={\bf b}_1\vee {\bf b}_2\vee {\bf b}_3$ has length $r_3$ (with some of the ${\bf b}_i$'s possibly equal to the empty set), ${\bf x}={\bf x}_1\vee {\bf x}_2\vee {\bf x}_3\vee {\bf x}_4$ and ${\bf t}_1\vee {\bf t}_2\vee{\bf t}_3\vee {\bf t}_4 = \sigma({\bf t})$ for some permutation $\sigma$. Squaring $F$ and integrating, we claim that there exist integers $s_1,s_2,s_3,s_4, s_5 \in \{1,...,q-1\}$ such that
\[
\|F\|^2_{\HH^{3q-2r_1-2r_2-2r_3}} = \int_{Z^{B}}\prod_{i=1}^5 (f\otimes_{s_i}f)_{\sigma_i}({\bf z}_{b_i})\mu(dz_1)\cdot\cdot\cdot\mu(dz_B),
\]
where $B=5q-2(s_1+ s_2+s_3+s_4+s_5)$, $\sigma_i$, $i=1,2,3,4, 5$, is a permutation of $[2q-2s_i]$, and the sets $b_1,b_2,b_3,b_4, b_5$ verify property {\bf ($\beta$)}, as defined at the beginning of this section. We have to consider separately different cases.

\noindent $(a)$: the length of ${\bf x}_4$ is not $0$. We can then consider separately the three first factors, for which the same expressions as in the proof of \eqref{E0} are available, and the two last ones, which give rise to  $s_4=s_5$ equal to  the length of ${\bf x}_4$.

\noindent $(b)$: the length of ${\bf x}_4$ is $0$ and the length of ${\bf t}_4$ is not $0$. Then  we consider separately the four  factors which are distinct from the fourth one and proceed as in the proof of \eqref{E1} for them, while the fourth one gives rise to $f\otimes_{\tau}f  $, with $\tau$ equal to the length of ${\bf t}_4$.

\noindent $(c)$: the lengths of ${\bf x}_4$  and  ${\bf t}_4$ are  $0$, but the length of ${\bf x}_3$ is not $0$. We then separate the five
factors into two groups, one with $f({\bf x_3},{\bf a}_1,{\bf a}_2,{\bf b}_3,{\bf t}_3)$ and $f({\bf x}_1,{\bf x}_2,{\bf x}_3)$, the other one with the three other factors. The first group gives rise to factors $f\otimes_{\tau}f  $, with $\tau$ equal to the length of ${\bf x}_3$, while the second group can be treated as in the proof of \eqref{E0}.

\noindent $(d)$: the lengths of ${\bf x}_3$, ${\bf x}_4$  and  ${\bf t}_4$ are  $0$, but the length of ${\bf t}_3$ is not $0$. We then consider separately the  factor  $f({\bf a}_1,{\bf a}_2,{\bf b}_3,{\bf t}_3)$, which gives rise to  a factor $f\otimes_{\tau}f  $, with $\tau$ equal to the length of ${\bf t}_3$. The four other factors can be treated as in the proof of \eqref{E1}.

\noindent $(e)$: the lengths of ${\bf x}_3$, ${\bf x}_4$, ${\bf t}_3$  and  ${\bf t}_4$ are  $0$. Remark that ${\bf x}_1$, ${\bf x}_2$ and
${\bf b}_3$ are non empty and, without loss of generality we can assume that ${\bf a}_2$ is non empty.  As before, we can conclude by
separating the five factors into two groups:
for the first one we take the first factor and $f({\bf x}_1,{\bf x}_2)$ whereas, for the second one, we choose the three remaining factors.

The desired conclusion (that is, (\ref{gogoal})) follows once again from Lemma \ref{l:CS}.

\fin

\section{Proof of the main results}\label{s:proofs}

\subsection{Proof of Theorem \ref{t:mainupper}}

The assumption $E[F^2]=1$
implies that $K:=1+E[|F|]\leq 2$.
The proof follows then immediately from (\ref{PBB2}) and (\ref{E1}).

\subsection{Proof of Theorem \ref{t:mainlower}}

Since $E[F_n]=0$ and $E[F_n^2]=1$, one has that $\kappa_4(F_n)=E[F_n^4]-3>0$.
Moreover, because $F_n\,\overset{{\rm Law}}{\rightarrow}\,N\sim\mathscr{N}(0,1)$ by assumption and due again to the hypercontractivity
of chaotic random variables, we have that $\kappa_4(F_n)=E[F_n^4]-E[N^4]\to 0$ as $n\to\infty$.
In the forthcoming proof we will need the following lemma.

\begin{lemme}\label{aline}
There exists $g,h\in \mathcal{U}\cap\mathcal{C}^\infty$ with bounded derivatives of all orders (except possibly the first one)
such that $E[f''_g(N)]\neq 0$, $E[f'''_g(N)]=0$, $E[f''_h(N)]=0$ and $E[f'''_h(N)]\neq 0$. 
\end{lemme}
{\it Proof}.
Let $H_p$, $p\geq 1$, denote the sequence of Hermite polynomials. Using the well-known formula
\[
\phi(x)=\sum_{p=0}^\infty \frac{1}{p!}E[\phi^{(p)}(N)]H_p(x),\quad x\in\R\,\mbox{ a.e.},\quad N\sim\mathcal{N}(0,1),
\]
valid for $\phi\in\mathcal{C}^\infty$ whose derivatives are all square integrable,
it is readily checked that, for almost all $x\in\R$,
\begin{eqnarray}\notag
\frac{\sqrt{e}}{1+\sqrt{e}}\cos x=\sum_{q=0}^\infty \frac{(-1)^q}{(2q)!(1+\sqrt{e})}H_{2q}(x)\quad\mbox{and}
\quad
\sin x=\sum_{q=0}^\infty \frac{(-1)^q}{(2q+1)!\sqrt{e}}H_{2q+1}(x).\\
\label{expan}
\end{eqnarray}
On the other hand, by applying several integration by parts,
we can write, for  $h\in\mathcal{U}$,
\begin{eqnarray}
E[f''_h(N)]&=&\int_{-\infty}^{+\infty}f''_h(x)\frac{e^{-x^2/2}}{\sqrt{2\pi}}dx=
\int_{-\infty}^{+\infty}f_h(x)(x^2-1)\frac{e^{-x^2/2}}{\sqrt{2\pi}}dx \notag\\
&=&\int_{-\infty}^{+\infty}dx\,
H_{2}(x)\int_{-\infty}^x dy\big(h(y)-E[h(N)]\big)\frac{e^{-y^2/2}}{\sqrt{2\pi}}\notag\\
&=&-\frac{1}{3}\int_{-\infty}^{+\infty} H_{3}(x)\big(h(x)-E[h(N)]\big)\frac{e^{-x^2/2}}{\sqrt{2\pi}}dx
=-\frac13\,E[h(N)H_3(N)].\notag\\
\label{ide}
\end{eqnarray}
Similarly, we can prove that, for all $h\in\mathcal{U}$,
\begin{equation}\label{ide2}
E[f'''_h(N)]=-\frac14\,E[h(N)H_4(N)].
\end{equation}
Now, let us consider $g(x)=\sin x$. Using (\ref{expan}) and then (\ref{ide})-(\ref{ide2}), we get
that $E[f''_g(N)]=\frac1{3\sqrt{e}}\neq 0$ and $E[f'''_g(N)]=0$. Moreover, $g$ belongs to $\mathcal{U}$ because $|g''(x)|=|\sin x|\leq 1$,
and has bounded derivatives.
Similarly, consider $h(x)=\frac{1}{1+\sqrt{e}}\big(\sqrt{e}\,\cos x -1+\frac12\,H_2(x)\big)$. 
Using once again (\ref{expan}) and then (\ref{ide})-(\ref{ide2}), we get this time
that $E[f''_h(N)]=0$ whereas $E[f'''_h(N)]=-\frac1{4+4\sqrt{e}}\neq 0$, 
also with $|h''(x)|=\left|\frac{1}{1+\sqrt{e}}\big(1-\sqrt{e}\,\cos x\big)\right|\leq 1$ so that
$h\in\mathcal{U}$.
\fin

\medskip

\noindent{\it Proof of Theorem \ref{t:mainlower}}. Recall that $\kappa_3(F_n)=E[F_n^3]$ and $\kappa_4(F_n)=E[F_n^4]-3$, and 
let $g,h\in\mathcal{U}$ be as in the statement of Lemma \ref{aline}. From Corollary \ref{c:est1} and Proposition \ref{prop-gamma}, we deduce
that
\begin{eqnarray*}
&&\left|
E[g(N)]-E[g(F_n)]-\frac12\,E[f''_g(N)]\kappa_3(F_n)
\right|\\
&\leq&\frac12\big|\kappa_3(F_n)\big|\left|E[f''_g(F_n)]-E[f''_g(N)]\right|+\frac16\left|E[f'''_g(F_n)]\right|\kappa_4(F_n)
+2c_4\|g'''\|_\infty\,\kappa_4(F_n)^{5/4}.
\end{eqnarray*}
Set
\[
c=\frac13\,\min\left\{ \frac12 |E[f''_g(N)]|, \frac16|E[f'''_h(N)]|\right\}.
\]
As $n\to\infty$, we have $E[f''_g(F_n)]- E[f''_g(N)]\to 0$, $E[f'''_g(F_n)]\to E[f'''_g(N)]=0$, and $\kappa_4(F_n)\to 0$.
Therefore, for $n$ large enough we have that
\[
d(F_n,N)\geq \big|E[g(N)]-E[g(F_n)]\big|\geq 3c|\kappa_3(F_n)| - \frac{c}{2}\max\{|\kappa_3(F_n)|,\kappa_4(F_n)\}.
\]
Similarly, we have 
\begin{eqnarray*}
&&\left|
E[h(N)]-E[h(F_n)]-\frac16\,E[f'''_h(N)]\kappa_4(F_n)
\right|\\
&\leq&\frac12\big|\kappa_3(F_n)\big|\left|E[f''_h(F_n)]\right|+\frac16\left|E[f'''_h(F_n)]-E[f'''_h(N)]\right|\kappa_4(F_n)
+2c_4\|h'''\|_\infty\,\kappa_4(F_n)^{5/4},
\end{eqnarray*}
from which we deduce, again for $n$ large enough, that
\[
d(F_n,N)\geq \big|E[h(N)]-E[h(F_n)]\big|\geq 3c\,\kappa_4(F_n) - \frac{c}{2}\max\{|\kappa_3(F_n)|,\kappa_4(F_n)\}.
\]
Finally taking the mean of the two previous upper bounds for $d(F_n,N)$ yields
\[
d(F_n,N)\geq \frac{3c}{2}\big(|\kappa_3(F_n)|+\kappa_4(F_n)\big) - \frac{c}{2}\max\{|\kappa_3(F_n)|,\kappa_4(F_n)\}
\geq c\,\max\{|\kappa_3(F_n)|,\kappa_4(F_n)\}.
\]
The proof is concluded.\fin

\section{Application: estimates in the Breuer-Major CLT}\label{s:BM}

In this final section, we determine optimal rates of convergence
associated with the well-known {\it Breuer-Major CLT} for Gaussian-subordinated random sequences -- see \cite{breuermajor}
for the original paper, or \cite{NPP} for a more modern reference. 
In order to be able to directly apply our previous results, we focus on sequences that can be represented as partial sums
of Hermite polynomials.

\subsection{General framework}\label{s:gf}

Consider a centered stationary Gaussian sequence $(X_k)_{k\in \mathbb{Z}}$
with unit variance and covariance function given by $E[X_{k}X_{l}] = \rho(k-l)$, $k,l\in\mathbb{Z}$.
Fix an integer $q\geq 2$, and set
\[
V_n = \frac{1}{\sqrt{n}}\sum_{k=0}^{n-1} H_q(X_k), \quad n\geq 1.
\]
Here, $H_q$ stands for the $q$th Hermite polynomial defined by (\ref{hq}).
Let also $N\sim\mathscr{N}(0,1)$,
and define $v_n :=E[V_n^2]$, $n\geq 1$. Finally, set
\[
F_n = \frac{V_n}{\sqrt{v_n}} = \frac{1}{\sqrt{n\,v_n}}\sum_{k=0}^{n-1} H_q(X_k).
\]
Without loss of generality, we may
assume that $X_k = X(h_k)$, where $X = \{X(h) : h\in \HH\}$ is some
isonormal Gaussian process and $\langle h_k,h_l\rangle_\HH=\rho(k-l)$ for every $k,l\in \mathbb{Z}$.
We then have
\[
F_n=I_q(f_n),\quad \mbox{with }
f_n=\frac{1}{\sqrt{nv_n}}\sum_{k=0}^{n-1} h_k^{\otimes q}.
\]

For $p \geq 1$, we introduce the Banach space $\ell^p(\Z)$ of $p$-summable sequences
equipped with the norm $\|u\|_{\ell^p} = \left(\sum_{k\in\Z}|u(k)|^p\right)^{1/p}$.
In what follows, we shall assume that $\rho$ belongs to $\ell^q(\Z)$.
Under this assumption, the celebrated Breuer-Major CLT (see \cite{breuermajor}, as well \cite[Chapter 7]{np-book}) asserts that 
\[
F_n\overset{\rm Law}{\to}\mathscr{N}(0,1)\quad\mbox{ as $n\to\infty$.}
\]

\begin{rem}{\rm 
One has that \[
v_n\to q!\sum_{k\in \Z}\rho(k)^q>0
\quad\mbox{as $n\to\infty$}.
\] 
It follows that, in the subsequent discussion, the role of the sequence $v_n$, $n\geq 1$, will be immaterial as far as rates of convergence are concerned.}
\end{rem}

\subsection{Explicit formulas for the third and fourth cumulants}

Let us compute,  in terms of $\rho$,  explicit expressions for the third and fourth cumulants of $F_n$. 
According to Proposition \ref{thm-pasmal}
(and using the notation introduced in Section \ref{s:steincum}), one has
\begin{equation}\label{Gamma1}
\Gamma_1(F_n)=q\sum_{r=1}^{q}(r-1)!\binom{q-1}{r-1}^2I_{2q-2r}
\left(f_n\widetilde{\otimes}_r f_n\right).
\end{equation}
Since $\kappa_3(F_n)=2E[F_n\Gamma_1(F_n)]$ by (\ref{azertiop}), we deduce
the following expression of $\kappa_3(F_n)$ in terms of the sequence $f_n$, $n\geq 1$:
\begin{equation}\label{k3}
\kappa_3(F_n)=\left\{\begin{array}{ll}2qq!(q/2-1)!\binom{q-1}{q/2-1}^2
\langle f_n,f_n\widetilde{\otimes}_{q/2}f_n\rangle_{\HH^{\otimes q}} &\mbox{ for even  } q\\
0 &\mbox{ for odd  } q\end{array}\right..
\end{equation}
Hence, when $q$ is even (observe that
$\langle f_n,f_n\widetilde{\otimes}_{q/2}f_n\rangle_{\HH^{\otimes q}}=\langle f_n,f_n\otimes_{q/2}f_n\rangle_{\HH^{\otimes q}}$
because $f_n$ is symmetric), we have
\begin{equation}\label{k3rho}
   \kappa_3(F_n)
=\frac{d_3(q)}{v_n^{3/2}\,n\sqrt{n}}\sum_{j,k,l=0}^{n-1}\rho(k-l)^{\frac q2}\rho(k-j)^{\frac q2}\rho(l-j)^{\frac q2},
\end{equation}
with $d_3(q)=2qq!(q/2-1)!\binom{q-1}{q/2-1}^2$. It will be often useful to transform the previous expression as follows: we have
\begin{eqnarray}
   \kappa_3(F_n) &=&
\frac{d_3(q)}{v_n^{3/2}\,n\sqrt{n}}\sum_{j=0}^{n-1}\sum_{k,l=-j}^{n-1-j}\rho(k-l)^{q/2} \rho(k)^{q/2} \rho(l)^{q/2}\nonumber \\
&=&\frac{d_3(q)}{v_n^{3/2}\sqrt{n}}
\sum_{k,l\in\Z}\eta_n(k,l)\rho(k-l)^{q/2} \rho(k)^{q/2} \rho(l)^{q/2}\label{K3:1},
\end{eqnarray}
where 
\begin{equation}\label{eta}
\eta_n(k,l)=\left (1-\frac{\max(k,l)_+}n + \frac{\min(k,l)_-}n\right)\mathbf{1}_{\{|k|<n,|l|<n\}}.
\end{equation}
\begin{rems}
{\rm 
\begin{enumerate} 
\item When $q$ is even, one has $\kappa_3(F_n)>0$ for all $n$; indeed,
\[
\sum_{j,k,l=0}^{n-1}\rho(k-l)^{\frac q2}\rho(k-j)^{\frac q2}\rho(l-j)^{\frac q2}=\frac{1}{(q/2)!}
E\left[\sum_{j=0}^{n-1}\left(\sum_{k=0}^{n-1}H_{q/2}(X_k)\rho(k-j)^{q/2}\right)^2\right].
\]
\item When $q=2$, one can even prove that $\Gamma_2(F_n)>0$ for all $n$ (recall that $\kappa_3(F_n)=2\,E[\Gamma_2(F_n)]$). 
\end{enumerate}
}
\end{rems}

Now, let have a look at the fourth cumulant.
Recall from \cite{nunugio} that
\begin{equation}\label{k4}
\kappa_4(F_n)=
\sum_{r=1}^{q-1} q!^2\binom{q}{r}^2\left\{ 
\|f_n\otimes_{r} f_n\|^2_{\HH^{\otimes 2(q-r)}}
+\binom{2q-2r}{q-r}\|f_n\widetilde{\otimes}_{r} f_n\|^2_{\HH^{\otimes 2(q-r)}}
\right\}.
\end{equation}
For the non-symmetrized contractions of (\ref{k4}), the link with $\rho$ is easily obtained; indeed, for any $r=1,\ldots,q-1$, we have:
\begin{equation}\label{convolution}
\|f_n{\otimes}_r f_n\|^2_{\HH^{\otimes 2(q-r)}}=\frac1{v_n^2\,n^2}\sum_{i,j,k,l=0}^{n-1}\rho(k-l)^r\rho(i-j)^r\rho(k-i)^{q-r}\rho(l-j)^{q-r}.
\end{equation}
On the other hand, we will actually face no problem due to symmetrized contractions.
Indeed, we may forget them when deriving lower bounds, whereas
we can use the inequality 
\begin{equation}\label{strat4}
\|f_n\widetilde{\otimes}_{r} f_n\|^2_{\HH^{\otimes 2(q-r)}}\leq
\|f_n\otimes_{r} f_n\|^2_{\HH^{\otimes 2(q-r)}}
\end{equation}
when dealing with upper bounds.
 
\subsection{Estimates for the third and fourth cumulants}

We start with the following result about the asymptotic behavior of the third cumulant of $F_n$.
Of course, by virtue of (\ref{k3}), only the case where $q$ is even must be considered.
\begin{prop}\label{k3:estimate} Assume that $q\geq 2$ is even. Then,
\[
{\kappa_3(F_n)}\leq  \frac{d_3(q)}{v_n^{3/2}\sqrt{n}}\left(\sum_{|k|<n}|\rho(k)|^{3q/4}\right)^2.
\]
Moreover, if $\rho\in \ell^{\frac{3q}{4}}(\zZ)$, then 
\begin{equation}\label{k3q2}
\sqrt{n}\,\kappa_3(F_n)\to \frac{d_3(q)}{q!^{3/2}}\sqrt{2\pi} \frac{\int_{\T}g_{_{q/2}}(t)^3dt}{\left(\int_{\T}g_{_{q/2}}(t)^2dt\right)^{3/2}}\quad\mbox{as $n\to\infty$},
\end{equation}
where
$\displaystyle g_{_{q/2}}(t):=\sum_{l\in\mathbb{Z}}\rho(k)^{q/2}e^{ikt}$ is almost everywhere positive on the torus $\T = \R \backslash (2\pi\Z) $.
\end{prop}
{\it Proof.} 
Recall the identity (\ref{K3:1}) and, for any $k\in\mathbb{Z}$ and any $n\geq 1$, set 
$\rho_n(k)=\rho(k){\bf 1}_{\{|k|< n\}}$ and $|\rho_n|(k)=|\rho_n(k)|$.
Since $0\leq \eta_n(k,l)\leq 1$ for all $n\geq 1$ and all $k,l\in\Z$,
we deduce that 
\[
\kappa_3(F_n)
\leq
\frac{d_3(q)}{v_n^{3/2}\sqrt{n}}
\sum_{l\in\mathbb{Z}}\big(|\rho_n|^{q/2}*|\rho_n|^{q/2}\big)(l)\big(|\rho_n|^{q/2}\big)(l).
\]
At this stage, let us recall the (well-known) Young inequality: if $s,p,p'\in[1,\infty]$ are such that $\frac1p+\frac1{p'}=1+\frac1s$, then
\begin{equation}\label{young}
\|u*v\|_{\ell^s}\leq \|u\|_{\ell^p}
\|v\|_{\ell^{p'}}.
\end{equation}
Hence, using first H\"{o}lder inequality and then Young inequality yield
\begin{eqnarray*}
\sum_{l\in\mathbb{Z}}\big(|\rho_n|^{q/2}*|\rho_n|^{q/2}\big)(l)\big(|\rho_n|^{q/2}\big)(l)
\leq \||\rho_n|^{q/2}*|\rho_n|^{q/2}\|_{\ell^3}\|\rho_n^{q/2}\|_{\ell^\frac{3}{2}}
\leq \|\rho_n^{q/2}\|_{\ell^\frac{3}{2}}^3,
\end{eqnarray*}
which proves the first statement of the proposition. 

Let us now further assume that $\rho\in \ell^{\frac{3q}{4}}(\zZ)$. It implies that 
$\rho\in \ell^{{q}}(\zZ)$ or, equivalently, that $\rho^{q/2}$ belongs to $\ell^{2}(\zZ)$.
Thus, the function $g_{_{q/2}}$ is well defined in $L^2(\T)$, as being the Fourier series with coefficients $\rho^{q/2}$.
In particular, we deduce from Bessel-Parseval equality that 
\[
v_n\to q!\sum_{k\in\Z}\rho(k)^q=\frac{q!}{2\pi}\int_{\T}g_{_{q/2}}(t)^2dt\quad\mbox{as $n\to\infty$}.
\]
We also have that $\rho^{q/2}*\rho^{q/2}$ belongs to $\ell^{2}(\zZ)$,
that is, that $\sum_{l\in\mathbb{Z}}\rho^{q/2}*\rho^{q/2}(l)\rho(l)^{q/2}$ is an absolutely convergent series whose value is given by 
$\frac{1}{2\pi}\int_{\T}g_{_{q/2}}(t)^3dt$ (Bessel-Parseval equality).  Then, using (\ref{k3rho})-\eqref{K3:1} and dominated convergence, we get
(\ref{k3q2}).
Finally, $\rho^{q/2}$ being a covariance sequence as well (those of the stationary sequence $\frac{1}{\sqrt{(q/2)!}}H_{q/2}(X_k)$),
its spectral density $g_{_{q/2}}$ is positive almost everywhere. (See e.g. \cite{doob}.)
\qed

\medskip

The situation for the fourth cumulant turns out to be not so easy, except when $q=2$. 

\begin{prop}\label{k4:estimate}
There exists $C>0$ such that, for all $n\geq 1$,
\begin{equation}\label{k4prop}
\kappa_4(F_n)\leq \left\{\begin{array}{ll} \frac C {v_n^2 n}\left(\sum_{|k|<n}|\rho(k)|^{2q/3}\right)^3 &\mbox{ if } q\leq 3\\
\frac{C}{v_n^2 n}\left(\sum_{k=-n+1}^{n-1}|\rho(k)|^{q-1}\right)^2\left(\sum_{k=-n+1}^{n-1}|\rho(k)|^{2}\right) &\mbox{ if } q\geq 3\end{array}\right..
\end{equation}
If $q=2$ and $\rho\in\ell^{4/3}(\Z)$, then
\begin{equation}\label{k4q2}
n\,\kappa_4(F_n)\to 24\pi  \frac{\int_{\T}g_{_1}(t)^4dt}{\left(\int_{\T}g_{_1}(t)^2dt\right)^2}\quad\mbox{as $n\to\infty$},
\end{equation}
where $g_1(t):=\sum_{l\in\Z}\rho(k)e^{ikt}$ is almost everywhere positive on the torus $\mathbb{T}$.
If $q\geq 3$ and $\rho\in\ell^{2}(\Z)$, then $\liminf_{n\to\infty}n\kappa_4(F_n)>0$.
\end{prop}
{\it Proof. } 
Thanks to (\ref{k4}), (\ref{convolution}) and \eqref{strat4}, we `only' have to estimate, for any $1\leq r\leq q-1$, 
\[
I(r)=\frac{1}{n}\sum_{i,j,k,l=0}^{n-1}|\rho(k-l)|^r|\rho(i-j)|^r|\rho(k-i)|^{q-r}|\rho(l-j)|^{q-r}.
\]
We immediately see that
\begin{eqnarray*}
I(r)\leq \frac{1}{n}\sum_{k, j=0}^{n-1}\left(|\rho_{n}|^r*|\rho_{n}|^{q-r}\right)^2(k-j)
 \leq  \sum_{j=-n}^n\left(|\rho_{n}|^r*|\rho_{n}|^{q-r}\right)^2(j).
\end{eqnarray*}
Let us first assume that $q=2$, so that $r=1$ necessarily. In this case,
\begin{equation}\label{i1}
I(1)\leq \||\rho_{n}|*|\rho_{n}|\|_{\ell^2}^2\leq \|\rho_{n}\|_{\ell^{\frac 43}}^4,
\end{equation}
where we have used Young inequality (\ref{young}) to get the last inequality. This proves (\ref{k4prop}) when $q=2$.
Assume now that $q\geq 3$. By symmetry, we may and will assume that $r\leq q/2$. 
Young inequality (\ref{young}) yields that
\begin{equation}\label{ir}
I(r)\leq \||\rho_{n}^r|*|\rho_{n}|^{q-r}\|^2_{\ell^2}\leq \|\rho_{n}^r\|_{\ell^2}^2\|\rho_{n}^{q-r}\|_{\ell^1}^2.
\end{equation}
This shows the desired bound when $r=1$. For the other values of $r$ (if any), we can make use of the log-convexity of the $\ell^p$ norms. 
More precisely, for $\alpha, \beta$ such that
$2r= 2(1-\alpha)+(q-1)\alpha$ and $q-r=2(1-\beta)+(q-1)\beta$, recall that
\begin{eqnarray*}
 \sum_{j\in\Z} |\rho_{n}(j)|^{2r} \leq  \|\rho_{n}\|_{\ell^2}^{2(1-\alpha)} \|\rho_{n}\|_{\ell^{q-1}}^{(q-1)\alpha}\quad\mbox{and}\quad
  \sum_{j\in\Z} |\rho_{n}(j)|^{q-r} \leq \|\rho_{n}\|_{\ell^2}^{2(1-\beta)} \|\rho_{n}\|_{\ell^{q-1}}^{(q-1)\beta}.
\end{eqnarray*}
We then conclude that (\ref{k4prop}) holds true for any $q\geq 3$ as well. 
To finish the proof of Proposition \ref{k4:estimate}, let us compute the limit of
$n\|f_n {\otimes}_{1} f_n\|^2_{\HH^{\otimes 2(q-1)}}$ under the additional assumption that
$\rho\in\ell^{4/3}(\Z)$ if $q=2$ and $\rho\in\ell^2(\Z)$ if $q\geq 3$. We can write
 \begin{eqnarray}
  n\|f_n {\otimes}_{1} f_n\|^2_{\HH^{\otimes 2(q-1)}}&=&\frac1{v_n^2\,n}\sum_{i,j,k,l=0}^{n-1}\rho(k-l)^{q-1}\rho(i-j)^{q-1}
\rho(k-i)\rho(l-j)\notag\\
   &=&\frac{1}{v_n^2}
\sum_{j,k,l\in\mathbb{Z}}
\eta_n(j, k,l)\rho(l-k)\rho(j)\rho(k)^{q-1}\rho(l-j)^{q-1}\label{k4k4k4},
\end{eqnarray}
where 
$\eta_n(j,k,l)=\left (1-\frac{\max(j,k,l)_+}n + \frac{\min(j,k,l)_-}n\right)
{\bf 1}_{\{|j|<n,|k|<n,|l|<n\}}$ 
is bounded by $1$ and tends to $1$ as $n\to\infty$. 
We know from (\ref{i1}) if $q=2$ and from (\ref{ir}) if $q\geq 3$ that 
the series 
\[
\sum_{j,k,l\in\mathbb{Z}}\rho(l-k)\rho(j)\rho(k)^{q-1}\rho(l-j)^{q-1}
\] 
is absolutely convergent 
under our assumption on $\rho$.
By dominated convergence, we get that
\[
  n\|f_n {\otimes}_{1} f_n\|^2_{\HH^{\otimes 2(q-1)}}\to 
\sum_{l\in\mathbb{Z}}\left(\rho*\rho^{q-1}\right)^2(l)=\|\rho*\rho^{q-1}\|_{\ell^2}^2\quad\mbox{as $n\to\infty$},
\]
implying in turn that $\liminf_{n\to\infty} n \kappa_4(F_n)>0$ as expected. 
Finally, when $q=2$ we have $v_n\to\displaystyle \frac{1}{\pi}\int_{\T}g_{1}(t)^2dt$, 
so that (\ref{k4q2}) holds true thanks to \eqref{k4}.

\qed

\subsection{Breuer-Major CLT}

The following result is the so-called Breuer-Major theorem (it is called this way in honor of the seminal paper
\cite{breuermajor}). For sake of completeness, we provide here a modern proof, that relies on the bounds 
(\ref{k4prop}) for $\kappa_4(F_n)$ and Theorem \ref{e:npupper}.

\begin{thm}[Breuer-Major]
Assume $\rho\in\ell^q(\mathbb{Z})$. Then $F_n$ converges to $\mathscr{N}(0,1)$ in total variation as $n\to\infty$.
\end{thm}
\medskip

\noindent
{\it Proof}.
By virtue of Theorem \ref{e:npupper}, it is (surprisingly) enough to prove that $\kappa_4(F_n)$ tends to $0$
as $n\to\infty$.
Let us first assume that $q=2$. Then, using H\"{o}lder inequality for $k>M$, we have, for any $n>M$,
\[
\frac1{{n}} \left(\sum_{k=-n+1}^{n-1}|\rho(k)|^{\frac{4}{3}}\right)^3 \leq
\frac{{4}}{{n}}\left(\sum_{|k|\leq M}|\rho(k)|^{\frac{4}{3}}\right)^3 +
4\left(\sum_{M<|k|<n}|\rho(k)|^{2}\right)^2.
\]
We can then conclude by a standard argument (choosing $M$ large enough). For $q\geq 3$, we proceed analogously, after having noticed that
\[
\frac1{n}\left(\sum_{k=-n+1}^{n-1}|\rho(k)|^{q-1}\right)^2\left(\sum_{k=-n+1}^{n-1}|\rho(k)|^{2}\right)
=\left(\frac{1}{n^{1/q}}\sum_{k=-n+1}^{n-1}|\rho(k)|^{q-1}\right)^2\left(\frac{1}{n^{1-2/q}}\sum_{k=-n+1}^{n-1}|\rho(k)|^{2}\right),
\]
with
\[
\frac{1}{n^{1/q}}\sum_{k=-n+1}^{n-1}|\rho(k)|^{q-1}\leq  2^{1/q}\left(\sum_{|k|<n}|\rho(k)|^{q}\right)^{1-1/q}\leq  2^{1/q}\|\rho\|_{\ell^q}^{q-1},
\]
and, for any $M<n$,
\[
\frac{1}{n^{1-2/q}}\sum_{k=-n+1}^{n-1}|\rho(k)|^{2}\leq \frac{1}{n^{1-2/q}}\sum_{|k|\leq M}|\rho(k)|^{2}+ 2^{1-2/q}\left(\sum_{M<|k|<n}|\rho(k)|^{q}\right)^{2/q}.
\]
Hence, by choosing $M$ large enough, we get that $\kappa_4(F_n)$ tends to $0$,
and the proof is concluded.
\qed

\subsection{The discrete-time fractional Brownian motion}

Let $B_H$ be a fractional Brownian motion with Hurst index $H\in(0,1)$.
We recall that $B_H=\{B_H(t)\}_{t\in\R}$ is a centered
Gaussian process with continuous paths such that
\[
E[B_H(t)B_H(s)] = \frac12\big(|t|^{2H}+|s|^{2H}-|t-s|^{2H}\big), \quad s,t\in\R.
\]
The process $B_H$ is self-similar with stationary increments, and we refer the reader to Nualart \cite{nualartbook}
and Samorodnitsky and Taqqu \cite{samotaqqu} for its main properties.
In this section, we offer fine estimates for $\kappa_3(F_n)$ and $\kappa_4(F_n)$ respectively, when the sequence $(X_k)_{k\in\Z}$
corresponds to increments of $B^H$, that is,
\begin{equation}\label{xk}
X_k=B_{H}(k+1)-B_H(k),\quad k\in\Z.
\end{equation}
The $X_k$'s are usually called `fractional Gaussian noise' in the literature, and are centered stationary Gaussian random variables with covariance
\[
\rho(k)=E[X_rX_{r+k}]=\frac{1}{2}\left(|k+1|^{2H}-2|k|^{2H}+|k-1|^{2H}\right),\quad r,k\in\Z.
\]
The covariance behaves asymptotically as
\begin{equation}\label{asymcov}
\rho(k)\sim H(2H-1)|k|^{2H-2}\quad\mbox{as $|k|\to\infty$}.
\end{equation}
In particular, when $H>\frac12$ we observe that, for $|k|$ large enough,
\begin{equation}\label{eses}
\rho(k)\geq H(H-\frac12)(1+|k|)^{2H-2}.
\end{equation}

In the following results, we let the notation introduced in Section \ref{s:gf} prevail, and we assume that the sequence $(X_k)_{k\in\Z}$ is given by
(\ref{xk}). We also use the following convention for non-negative sequences $(u_n)$ and $(v_n)$ (possibly depending on $q$ and/or $H$): we write $v_n\asymp u_n$
to indicate that $0<\liminf_{n\to\infty}v_n/u_n\leq \limsup_{n\to\infty}v_n/u_n<\infty$.

\begin{prop} \label{sharp:k3} 
Assume that $q\geq 2$ is even. We have:
\[
\kappa_3(F_n)\asymp
\left\{\begin{array}{lll}
n^{-\frac 12} \;&\mbox{if } &0<H<1-\frac2{3q} \\
n^{-\frac 12}\log^2 n \;&\mbox{if } &H=1-\frac2{3q} \\
n^{\frac32-3q+3qH}\;&\mbox{if }& 1-\frac2{3q}<H<1-\frac{1}{2q} \end{array}\right..
\]
\end{prop}
{\it Proof.} 
When $H<1-\frac2{3q}$, we have that $\rho \in \ell^{3q/4}(\Z)$, so that 
\[
0<\liminf_{n\to\infty}\sqrt{n}\,\kappa_3(F_n)=\limsup_{n\to\infty}\sqrt{n}\,\kappa_3(F_n)<\infty
\] 
by Proposition \ref{k3:estimate}.
For the other values of $H$, consider first the limsup. We have that 
\[
\sum_{|k|<n}|\rho(k)|^{3q/4}\leq C \left\{\begin{array}{ccl}
\log n \;&\mbox{if } &H=1-\frac2{3q}\\
n^{1-\frac{3q}2+\frac{3qH}2}\;&\mbox{if }& H >1-\frac2{3q} \end{array}\right.,
\]
from which we deduce the finiteness of the limsup by Proposition \ref{k3:estimate}
together
with the fact that $v_n\to q!\sum_{k\in\Z}\rho(k)^q\in(0,\infty)$ when $H<1-\frac{1}{2q}$.

Let us now focus on the liminf. 
Since $H>1-\frac{2}{3q}>\frac12$, the lower bound (\ref{eses}) holds for $|k|$ large enough.
In fact, by considering a decomposition of the type $\rho=\rho_+ + \tau$
with $\tau$ having compact support,
we can do as if the inequality (\ref{eses}) were valid for the small values of $|k|$ as well. Indeed, it suffices to
use the fact that $\kappa_3(F_n)$ comes from a trilinear form (see (\ref{K3:1})) as well as the inequality $|x|^{q/2}|y|^{q/2}\leq
\frac12\big(|x|^q+|y|^q\big)$, so to deduce that, when 
there is at least one term $\tau$, the asymptotic order of the corresponding contribution in $\kappa_3(F_n)$ is $O(n^{-\frac 12})$.
Recall the definition (\ref{eta}) of $\eta_n(k,l)$.
With the extra assumption that (\ref{eses}) holds true for all $k$, we get that 
\[
\sum_{k,l\in\Z}\eta_n(k,l)\rho(k-l)^{q/2} \rho(k)^{q/2} \rho(l)^{q/2}
\]
is bounded by below by a constant time
\[
\sum_{|l|\leq n/4}\sum_{|l|\leq |k|\leq 2|l|}(1+|k-l|)^{(H-1)q}(1+|k|)^{(H-1)q} (1+|l|)^{(H-1)q}.
\]
Note that, for $k,l\in\Z$ with $|l|\leq |k|\leq 2|l|$, we have
$1+|k-l|\leq 2(1+|k|)$ as well as $1+|k-l|\leq 3(1+|l|)$, so that
\[
(1+|k-l|)^{(H-1)q}\geq 6^{(H-1)q/2}(1+|k|)^{(H-1)q/2} (1+|l|)^{(H-1)q/2},
\]
which, by  summing first over $k$ and then over $l$, concludes the proof for the liminf.
\qed

\medskip

By reasoning similarly, we obtain an estimate for $\kappa_4(F_n)$. 
(In the following statement, the limsup is partially known from \cite{BBL}.)

\begin{prop} \label{sharp} 
For $q\in\{2,3\}$, we have
\begin{equation}\label{rate-quatre}
\kappa_4(F_n)\asymp
\left\{\begin{array}{lll}n^{-1}\;&\mbox{if } & 0<H  <1-\frac{3}{4q} \\
n^{-1}\log^3 n \;&\mbox{if } &H=1-\frac{3}{4q} \\
n^{4qH-4q+2}\;&\mbox{if }& 1-\frac{3}{4q}<H<1-\frac{1}{2q} \end{array}\right.\;.
\end{equation}
whereas, for $q>3$,
\begin{equation}\label{rate-qBFM}
\kappa_4(F_n)\asymp
\left\{\begin{array}{lll}
n^{-1}\;&\mbox{if }& 0<H <\frac34\\
n^{-1}\log (n)\;&\mbox{if }& H =\frac34 \\
n^{4H-4}\;&\mbox{if } &\frac34<H<1-\frac1{2q-2}\\
n^{4H-4}\log^2 n\;&\mbox{if } & H  = 1-\frac 1{2q-2}\\
n^{4qH-4q+2}\;&\mbox{if } &  1-\frac 1{2q-2}<H<1-\frac{1}{2q}
\end{array}\right.\;.
\end{equation}\end{prop}
{\it Proof.}
The proof of the limsup is straightforward by using Proposition \ref{k4:estimate}. 
To get the liminf result, we first consider $\rho=\rho_+ +\tau$ 
as in the proof of Proposition \ref{sharp:k3}, 
and verify that the contribution of all the terms containing at least one $\tau$ 
are of lower order. 
Here again, we can therefore do as if the inequality (\ref{eses}) were valid for all $k$.
This allows us to bound (\ref{k4k4k4}) by below by following the same line of reasoning as in the proof of Proposition \ref{sharp:k3}. 
We finally conclude thanks to (\ref{k4}).
Details are left to the reader.
\qed

\medskip

By comparing the asymptotic behaviors of $\kappa_3(F_n)$ and $\kappa_4(F_n)$, we observe the following non-expected fact.

\begin{cor}\label{cor-surprising} 
When $q\geq 6$, there exists 
a non-trivial range of values of $H$ for which $\kappa_4(F_n)$ tends less rapidly to $0$ than $\kappa_3(F_n)$.
\end{cor}
{\it Proof}.
The two functions that give the behavior of $\kappa_3(F_n)$ and $\kappa_4(F_n)$ are piecewise linear and concave. It is therefore
sufficient to consider the value 
$H=1-\frac2{3q}$. For this value, one has
$\kappa_4(F_n)\asymp n^{4H-4}\gg n^{-1/2}$ when $\frac 1{q-1}<\frac 4{3q}< \frac 14$.
\qed


\begin{thebibliography}{99}

\bibitem{BarbourPtrf1986}
A.D. Barbour (1986).
Asymptotic expansions based on smooth functions in the central limit theorem.
{\it Probab. Theory Rel. Fields} {\bf 72}(2), 289-303.


\bibitem{BBL}
H. Bierm\'e, A. Bonami and J. L\'eon (2011).
Central Limit Theorems and Quadratic Variations in terms of  Spectral Density.
{\it Electron. J. Probab.} {\bf 16}, 362-395.


\bibitem{breuermajor}
P. Breuer and P. Major (1983). Central limit theorems for non-linear functionals of Gaussian fields. {\it J. Mult. Anal.} {\bf 13}, 425-441.

\bibitem{ChenGoldShao} L. H. Y. Chen, L. Goldstein and Q.-M. Shao (2011). {\it Normal Approximation by Stein's Method}. Springer-Verlag, Berlin

\bibitem{daly}
F.A. Daly (2008).
Upper bounds for Stein-type operators.
{\it Electron. J. Probab.} {\bf 13}, paper 20, 566-587.

\bibitem{doob}
J.L. Doob (1953).
\it Stochastic Processes.
\rm Wiley.

\bibitem{Janson} S. Janson (1997). \textit{Gaussian Hilbert Spaces. }
Cambridge University Press, Cambridge.

\bibitem{MallBook} P. Malliavin (1997). \textit{Stochastic Analysis}.
Springer-Verlag, Berlin, Heidelberg, New York.

\bibitem{np-ptrf}
\rm I. Nourdin and G. Peccati (2009):
\rm Stein's method on Wiener chaos.
\it Probab. Theory Relat. Fields \rm {\bf 145}, no. 1, 75-118.

\bibitem{np-aop}
\rm I. Nourdin and G. Peccati (2009):
\rm Stein's method and exact Berry-Esseen asymptotics for functionals of Gaussian fields.
\it Ann. Probab. \rm {\bf 37}, no. 6, 2231-2261.

\bibitem{survey}
\rm I. Nourdin and G. Peccati (2010):
\rm Stein's method meets Malliavin calculus: a short survey with new estimates.
\rm In the volume: \it Recent Development in Stochastic Dynamics and Stochastic Analysis, \rm World Scientific, 207-236.

\bibitem{np-jfa}
\rm I. Nourdin and G. Peccati (2010):
\rm Cumulants on the Wiener space.
\it J. Funct. Anal. \rm {\bf 258}, 3775-3791.

\bibitem{np-book} I. Nourdin and G. Peccati (2011): {\it Normal Approximations with Malliavin Calculus. From Stein's Method to Universality}. Cambridge University Press, forthcoming book.

\bibitem{NPP}
I. Nourdin, G. Peccati and M. Podolskij (2011).
Quantitative Breuer-Major Theorems.
\it Stoch. Proc. Appl. \rm {\bf 121}, no. 4, 793-812. 


\bibitem{npr-aop}
I. Nourdin, G. Peccati and G. Reinert (2010). Invariance principles for homogeneous sums: universality of Gaussian Wiener chaos. {\it Ann. Probab.} {\bf 38}(5), 1947-1985.

\bibitem{nualartbook}
\rm D. Nualart (2006).
\it The Malliavin calculus and related topics of Probability and Its Applications.
\rm Springer-Verlag, Berlin, second edition.

\bibitem{nunugio}
\rm D. Nualart and G. Peccati (2005):
\rm Central limit theorems for sequences of multiple stochastic integrals.
\it Ann. Probab. \rm {\bf 33}, no. 1, 177-193.

\bibitem{PeTa}
G. Peccati and M.S. Taqqu (2010). {\it Wiener chaos: moments, cumulants and diagrams}. Springer-Verlag.

\bibitem{RotarOnBar} V. Rotar (2005). Stein's method, Edgeworth's expansions and a formula of Barbour.
In: {\it An Introduction to Stein's Method} (A.D. Barbour and L.H.Y. Chen, eds), Lecture Notes Series No.{\bf  4},
Institute for Mathematical Sciences, National University of Singapore, Singapore University Press and World Scientific 2005, 59-84.

\bibitem{samotaqqu}
G. Samorodnitsky and M.S. Taqqu (1994). {\it Stable Non-Gaussian Random Processes}. Chapman and Hall, New York.

\bibitem{Shigekawa} I. Shigekawa (1980). Derivatives of Wiener functionals and absolute continuity of induced measures. {\it J. Math. Kyoto Univ.} {\bf 20}(2), 263-289.

\end{thebibliography}
\end{document}